\documentclass[journal]{IEEEtran}
\usepackage{graphicx} 
\usepackage{subcaption}
\usepackage{tabularray}
\usepackage{cite}
\usepackage{amsmath, amssymb, amsfonts, enumerate, bbm, bm}
\usepackage[ruled]{algorithm2e}
\usepackage{mathrsfs}

\def\reals{\mathbb R}

\def\Tcal{\mathcal T}

\def\ta{t^\mathrm{A}}
\def\ts{t^\mathrm{S}}
\def\te{t^\mathrm{E}}
\def\td{t^\mathrm{D}}

\newtheorem{theorem}{Theorem}
\newtheorem{lemma}{Lemma}

\newtheorem{corollary}{Corollary}

\newtheorem{remark}{Remark}
\newtheorem{example}{Example}

\newtheorem{assumption}{Assumption}
\usepackage[usenames,dvipsnames,table]{xcolor}

\newcommand{\rev}[1]{{\color{black} #1}}
\newcommand{\orange}[1]{{\color{black} #1}}
\newcommand{\blue}[1]{{\color{black} #1}}

\newcommand{\brown}[1]{{\color{black} #1}}
\newcommand{\purple}[1]{{\color{black} #1}}
\newcommand{\gray}[1]{{\color{black}#1}}
\newcommand{\orangeV}[1]{{\color{black} #1}}

\usepackage[normalem]{ulem}

\usepackage{tikz}
\usetikzlibrary{automata}
\usetikzlibrary{positioning, arrows,calc}
\tikzset{
	->,  
	>=stealth, 
	shorten >=2pt, shorten <=2pt, 
	node distance=3cm, 
	every state/.style={draw=blue!55,very thick,fill=blue!20}, 
	initial text=$ $, 
 }

\usepackage{hyperref}
\hypersetup{
	colorlinks,
	linkcolor={red!50!black},
	citecolor={blue!50!black},
	urlcolor={blue!50!black}
}

\begin{document}
\title{Extended Version: Non-Preemptive Scheduling of Flexible Loads in Smart Grids via Convex Optimization}

\author{Mehdi Davoudi, Mingyu~Chen, %
        and~Junjie~Qin
\thanks{M. Davoudi and J. Qin are with the Elmore School of Electrical and Computer Engineering, Purdue University, West Lafayette, IN, US. Email:
        {\tt\small \{mdavoudi,jq\}@purdue.edu}.}%
\thanks{M. Chen is with the Department of Electrical and Computer Engineering, Boston University. Email:
        {\tt\small mingyuc@bu.edu}. }
\thanks{This work is supported in part by a grant from the Schweitzer Engineering Laboratories and in part by the U.S. National Science Foundation under Grant No. ECCS-2339803. The authors are grateful to invaluable comments and suggestions offered by Dr. Ellery Blood (Schweitzer Engineering Laboratories), Dr. Dionysios Aliprantis (Purdue) and Dr. Steven Pekarek (Purdue). }
}


\maketitle

\addtolength{\textfloatsep}{-2mm}
\begin{abstract}

This paper studies the scheduling of a large population of non-preemptive flexible electric loads, each of which has a flexible starting time but once started will follow a fixed load shape until completion. We first formulate the scheduling problem  as a mixed-integer convex program (MICP), then propose an efficient polynomial time relaxation-adjustment-rounding algorithm for solving the problem. The key novelty of the proposed method lies in its adjustment step, which uses a  graph-based algorithm to navigate within the set of optimal points of the convex relaxation while reducing the number of fractional entries in the solution. We establish mathematically that our algorithm yields solutions that are near optimal for a finite number of loads and with its sub-optimality independent of the number of loads. Consequently, the proposed method is asymptotically optimal in a per-load cost sense when the number of loads increases. Despite the gap between the MICP and its convex relaxation, we establish that the solution of the proposed algorithm can be decentralized by marginal prices of the convex relaxation. We also develop and analyze variants of the proposed algorithm for settings with uncertainty and with time-varying realistic load shapes. Finally, we numerically evaluate the proposed algorithm in a case study for the non-preemptive scheduling of electric vehicles charging loads. 

\end{abstract}
\section{Introduction}
 The ongoing trend of replacing conventional thermal generation units with renewable sources \purple{that are} characterized by \purple{their} inherently intermittent output, has reduced the controllability of the supply side of the power system \cite{lund2015review}. There are various solutions to address this challenge, among which a potentially cost-effective one resides on the demand side, where the adoption of communication, control, and computing technologies on the grid edge (e.g., building management systems, smart home technology, \purple{and managed charging solutions}) can unlock  substantial  demand-side flexibility \cite{morstyn2018designing}. However, leveraging this flexibility necessitates coordination \purple{of} diverse loads with various constraints and owner \purple{preferences} \cite{fruh2022coordinated}, which is challenging in practice.  

 There are different types of flexible loads, among which shiftable loads play a major role in providing demand-side flexibility \cite{tang2021flexibility}. Shiftable loads can be categorized into two groups. The first one is \emph{preemptive shiftable loads} in which not only their starting time is flexible but also they can be interrupted and have adjustable power consumption profiles. The scheduling of these loads may be effectively computed by solving convex optimization problems for which efficient polynomial-time algorithms exist. Consequently, a significant portion of prior research on demand-side management has concentrated on preemptive shiftable loads \cite{shao2023preemptive,tu2019equivalent,li2021coordinating,qayyum2015appliance}. 

 The second category encompasses \emph{non-preemptive shiftable loads} whose starting time is flexible but once started they should follow a certain power consumption pattern until completion. Non-preemptive shiftable loads are prevalent and can be found in various applications, including numerous industrial processes, and residential loads such as washing machines and dishwashers \cite{dewangan2023improved}. Moreover, many preemptive shiftable loads without advanced communication and control functionalities can be controlled as non-preemptive loads using low-cost networked switches. 
 For instance, flexible electric vehicle (EV) charging loads may be modeled as non-preemptive when the charger only offers on/off control.
 However, scheduling a large population of non-preemptive loads remains a challenging open problem with limited prior exploration. This challenge arises from the inherent combinatorial nature of the underlying optimization problem, which in general is NP-hard \cite{korte2011combinatorial}.

\subsection{Related literature}\label{reveiw}
The predominant body of existing literature in this domain formulates the scheduling of non-preemptive loads as  mixed-integer  programs and then either relies on commercial solvers that implement branch and bound algorithms \cite{van2018techno,antunes2022comprehensive,sou2011scheduling}, or uses global optimization methods \cite{javaid2017new,zhu2019optimal} to solve the  problems.  
Similar formulations also arise in the context of EV charging scheduling, when on/off or discrete charging rate constraints are modeled \cite{sun2016optimal, binetti2015scalable}.
 Despite the contributions made by these studies, their proposed methods may suffer from a high computational complexity and/or lack provable performance guarantees. 
 

Several prior studies aim to address this issue by developing polynomial time algorithms with performance guarantees.
 O'Brien and Rajagopal \cite{7066982} propose a greedy algorithm  to schedule non-preemptive loads and mathematically establish that the algorithm is $(1-1/e)$-optimal.
 Non-preemptive load scheduling is studied in the context of matching market design for distributed energy resources in \cite{qin2018automatic}, for which a fluid relaxation based method is proposed and performance guarantees are proved for linear cost functions. 
Gupta et al. \cite{7447118} investigate a continuous-time version of the non-preemptive load scheduling problem and analyze the performance of the Earliest Deadline First (EDF) heuristic focusing on the case with time-invariant cost functions and load shapes. However,
their performance guarantee relies on certain assumptions
involving endogenous parameters derived from the algorithm
itself, which cannot be verified prior to execution.
 A market consisting of different agents including non-preemptive loads is considered in \cite{dahlin2022scheduling}, where the competitive equilibrium is obtained using the relaxed version of the non-preemptive scheduling problem. 
Their optimality results rely on a probabilistic allocation that is justified by assuming each load represents a population of identical loads.
 For EV charging problems with discrete charging rates, Gan et al. \cite{gan2012stochastic} devise a stochastic distributed algorithm and establish a sub-optimality bound that increases with the number of EVs.

 
While these papers laid the foundation for developing provably near optimal efficient algorithm for the non-preemptive load scheduling problem, the known sub-optimality bounds all increase with the number of loads, which  may scale unfavorably for large problems. \textcolor{black} {See Table \ref{compare} for a more detailed comparison.} Therefore, it is unclear whether these algorithms are well suited for practical settings with a large number of loads/EVs.
 It remains an open problem to develop an efficient (i.e., polynomial-time) algorithm with sub-optimality independent of the number of loads \textcolor{black}{for general non-preemptive load scheduling problems with  time-varying cost functions and load shapes.} 

\textcolor{black}{It is worth noting that combinatorial scheduling problems have been studied in other domains, including variants of machine scheduling \cite{lenstra1990approximation,phillips1997task,im2014dynamic} and broadcast scheduling \cite{im2014new}. 
Although  algorithms developed for these setups may seem relevant, 
the time-varying cost functions and load shapes, which are important to be considered in our load scheduling problem,
makes it difficult to directly apply existing scheduling algorithms to our problem. 
In addition to the need for new algorithms, the generality associated with time-varying costs and load shapes also demands different techniques for algorithm analysis, as the classical competitive ratio analysis in existing scheduling papers is not suitable for cost minimization problems without assuming the cost is bounded away from zero.}

 \subsection{Organization and contributions}
 In response to this need, we investigate the setting where an aggregator schedules a large number of non-preemptive shiftable loads (referred to as jobs). 
 We first formulate the scheduling problem as a mixed-integer convex program (MICP)  (Section~\ref{sec:model}) and then propose a relaxation-adjustment-rounding solution approach that is constituted of three steps (Section~\ref{sec:scheduling}): (a) We form and solve the convex relaxation of the MICP by replacing integer variables with continuous ones. The solution that we obtain may have some fractional entries, violating integer constraints. (b) We introduce a \emph{lossless} adjustment step based on a graph-based polynomial-time algorithm to decrease the number of fractional entries in the obtained solution, yielding another \emph{optimal point} for the convex relaxation problem. (c) The remaining fractional entries in the solutions are then eliminated via rounding. 
We establish mathematically that the proposed algorithm results in a schedule that is admissible  for each job. Additionally, we derive a bound \purple{for} the sub-optimality of the solutions which is independent of the number of jobs.
\purple{As a result}, our algorithm achieves an asymptotically optimal per-job cost as the number of jobs increases. 
We then generalize the algorithm to incorporate a number of practical considerations (Section~\ref{sec:ext}), including uncertainty, decentralized implementation, and realistic time-varying load shapes. 
The algorithm's performance is compared against directly solving the MICP with a commercial solver and the greedy heuristic in numerical experiments for the non-preemptive EV charging scheduling problem (Section~\ref{sec:ne}).

We contribute to the literature in the following ways:
\begin{enumerate}[(a)]
	\item To the best of our knowledge, this paper is the \emph{first} to propose a polynomial-time algorithm for scheduling non-preemptive loads with a sub-optimality independent of the number of loads (and thus asymptotically optimal in per-job cost sense) for general convex costs. 
	\item In doing so, our paper, for the first time, justifies the  folklore that large scale non-preemptive load scheduling problems can be well approximated by convex relaxation via offering and analyzing a novel algorithm (which is different from the standard relaxation-rounding procedure). The necessity of our adjustment step is highlighted by proving that the worst-case sub-optimality of the standard relaxation-rounding procedure is unbounded. 
	\item Despite the gap between the MICP formulation and its convex relaxation, we establish a somewhat surprising self-scheduling property that the solution of the proposed algorithm can be decentralized with marginal prices of the convex relaxation.
	\item We also develop practical variants of the proposed 
	algorithm and their associated performance guarantees for the case with uncertainty and realistic loadshapes, paving the way for real-world implementation.
\end{enumerate}

\begin{table*}[h]
	\centering
	\caption{\textcolor{black}{Comparison of non-preemptive load scheduling algorithms in the literature}}
	\label{compare}
	\renewcommand{\arraystretch}{1.3} 
	\setlength{\tabcolsep}{5pt} 
	\begin{tabular}{|>{\centering\arraybackslash}p{1.5cm}|>{\centering\arraybackslash}p{3.2cm}|>{\centering\arraybackslash}p{4cm}|>{\centering\arraybackslash}p{6.8cm}|}
		\hline
		\textcolor{black}{\textbf{Reference}} & \textcolor{black}{\textbf{Objective function}} & \textcolor{black}{\textbf{Load shape}} & \textcolor{black}{\textbf{Performance guarantee}} \\ 
		\hline\hline
		\textcolor{black}{{\cite{7066982}}} & \textcolor{black}{Time-invariant convex} & \textcolor{black}{Time-varying} & \textcolor{black}{($1 - 1/e$)-approximation} \\ 
		\textcolor{black}{{\cite{qin2018automatic}}} & \textcolor{black}{Time-invariant linear} & \textcolor{black}{Time-varying} & \textcolor{black}{Asymptotically optimal} \\ 
		\textcolor{black}{{\cite{7447118}}} & \textcolor{black}{Time-invariant convex} & \textcolor{black}{Time-invariant} & \textcolor{black}{Sub-optimality bound independent of number of loads} \\ 
		\textcolor{black}{{\cite{dahlin2022scheduling}}} & \textcolor{black}{Time-varying linear} & \textcolor{black}{Time-invariant} & \textcolor{black}{Asymptotically optimal} \\ 
		\textcolor{black}{{\cite{gan2012stochastic}}} & \textcolor{black}{Time-varying convex} & \textcolor{black}{Time-invariant} & \textcolor{black}{Sub-optimality bound dependent on number of loads} \\ 
		\textcolor{black}{This paper} & \textcolor{black}{Time-varying convex} & \textcolor{black}{Time-invariant \& Time-varying} & \textcolor{black}{Sub-optimality bound independent of number of loads} \\ 
		\hline
	\end{tabular}
\end{table*}

 This paper generalizes our prior conference paper\cite{chen2022scheduling} by (a) developing and analyzing new variants of the relaxation-adjustment-rounding algorithm for the  cases with realistic load shapes and uncertainty,   (b) providing more detailed exposition and complete proofs, and (c) performing extensive numerical experiments to empirically evaluate the proposed algorithms.
 

\section{Model}\label{sec:model} 

Consider an aggregator whose role is to schedule 
a large collection of non-preemptive shiftable loads. \purple{Each element of this collection is referred to as \gray{a} \emph{job} and denoted by $j\in \mathcal{J}:=\{1,\dots, J\}$, in which $J$ is the total number of non-preemptive loads.} We work
with a finite horizon discrete-time model, and denote each
time period by $t\in \Tcal :=\{1, \dots, T\}$\orange{,} where $T$ is the total
number of time periods. 

\vspace{.05in}
\noindent{\bf Notations:} For any job-specific and time-varying variable $x_j(t)$ and time-varying variable $y(t)$, $j \in \mathcal J$ and $t\in \mathcal T$, let $\mathbf x \in \reals^{J\times T}$ denote\gray{s} the collection  $\{x_j(t)\}_{j\in \mathcal J, t \in\mathcal T}$, $\mathbf x_j \in \reals^{T}$ denote\gray{s} the collection $\{x_j(t)\}_{t\in \mathcal T}$, $\mathbf x(t) \in \reals^{J}$ denote\gray{s} the collection $\{x_j(t)\}_{j\in \mathcal J}$, and $\mathbf y \in \reals^T$ denote\gray{s} the collection $\{y(t)\}_{t\in \mathcal T}$. We use $\mathbf 1$ and $\mathbf e_i$ to denote the all-one vector and the  $i$-th elementary vector of appropriate dimension, respectively.  

\subsection{Non-preemptive jobs}\label{model:def:job}
\def\Jcal{\mathcal J}
\def\Tcals{\Tcal^\mathrm{S}}
Ahead of time $t=1$, the aggregator is given $J$ non-preemptive electric loads. 
\gray{For simplicity, we first consider the case without uncertainty and with rectangular load shapes.} \purple{Cases with uncertain parameters and general load shapes are considered in Sections~\ref{stochastic} and~\ref{sec:realload}, respectively}.
\purple{A rectangular job $j \in \Jcal$} is characterized by parameters $(p_j,\, d_j,\, \ta_j,\, \td_j)$, where $p_j \in [0,\, p^{\max}]$ is the power requirement for running the job for each time period, $d_j\in \{1,\, \dots, d^{\max}\}$ is the duration of the job, $\ta_j\in \Tcal$ is the earliest time period in which the job is \emph{available} to start, and $\td_j\in \Tcal$ is the latest time period in which the job must finish, i.e., the \emph{deadline} for the job. 
Here $p^{\max}$ and $d^{\max}$ are the maximum power consumption level and the maximum duration among all the jobs, respectively. 
We assume there is a feasible schedule within $\Tcal$ for each job $j$, i.e., 
$1 \le \ta_j\le \ta_j+d_j-1\le \td_j \le T$. \purple{The length of each time period is $\Delta t$}. For notational simplicity, we focus on the case $\Delta t=1$ hour throughout this paper, thus consumed energy \purple{(kWh)} within each time interval and the \purple{average} power consumption \purple{(kW)} \gray{have the same numerical value.} Other $\Delta t$ values can be accommodated with a minimal change of the proposed method and results.

\begin{figure}[h]
	\centering
	\includegraphics[width=0.45\textwidth]{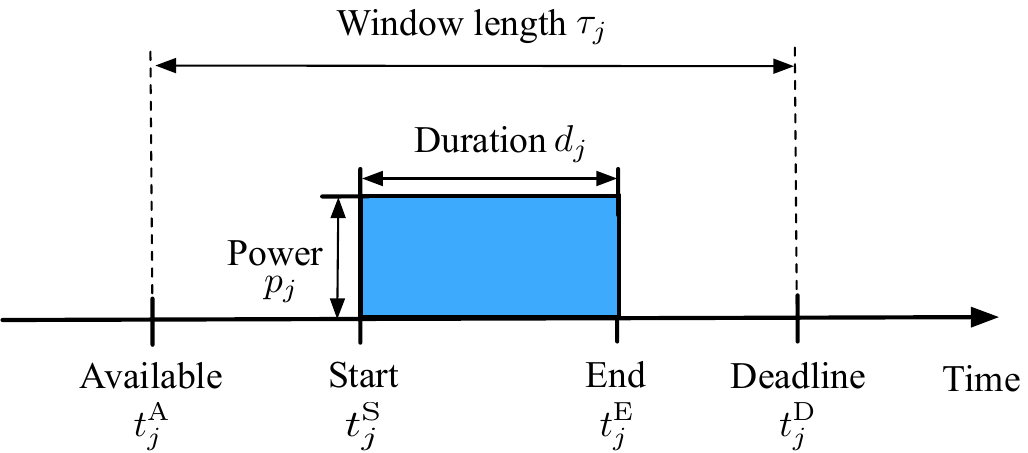}
	\caption{Illustration of a non-preemptive shiftable load (job)}\label{fig:jobw}
\end{figure}

We refer to the time window from $\ta_j$ to $\td_j$ as the \emph{availability window} of job $j$, and denote its length by 
\begin{equation}
\tau_j := \td_j +1 - \ta_j, \quad j\in \Jcal. \nonumber 
\end{equation}
A scheduling decision for job $j$ is deemed \emph{admissible} if the job is scheduled within its availability window, i.e., 
\begin{equation}
\ta_j \le \ts_j < \te_j \le \td_j, \quad j \in \Jcal,\nonumber
\end{equation}
where $\ts_j$ and $\te_j$ denote the time periods in which the job is scheduled to start and end, respectively.  
We also denote the \emph{set of admissible starting times} of job $j$ by 
\begin{equation}
\Tcals_j := \{\ta_j, \dots, \td_j-d_{\orange{j}}+1\}, \quad j \in \Jcal.\nonumber 
\end{equation}
Fig.~\ref{fig:jobw} provides an illustration of the parameters of a job.

Given a scheduling decision for job $j$, which can be represented by a starting time $\ts_j$, the power consumption profile of the job over $\mathcal T$, denoted by $\bm \ell_j \in \reals^T$,  is 
\begin{equation}
\ell_j(t) = \begin{cases}
	p_j,&  \mbox{if } \ts_j \le t \le \te_j, \\
	0,  &\mbox{otherwise,}
\end{cases}\nonumber
\end{equation}
where $\te_j= \ts_j+d_j-1$. 



\subsection{Cost model}
Given the job schedules, the aggregate power consumption profile of all the jobs at each time interval is 
\begin{equation}
L(t) = \sum_{j\in \Jcal} \ell_j(t), \quad t\in \Tcal.\nonumber 
\end{equation}
We model the cost of serving this aggregate load by
\begin{equation} \label{MILP_obj}
\Phi(\blue{\mathbf{L}}) \blue{:}= \sum_{t\in \Tcal} \phi_t(L(t)),
\end{equation}
where  the cost function for each period is $\phi_t(\cdot): \reals\mapsto \reals$, which is  convex  and $K$-Lipschitz continuous over its domain, with some finite Lipschitz coefficient $K>0$.   \rev{Depending on the actual setup (e.g., the load is served by purchasing electricity from wholesale electricity market or by generating with onsite generators), the cost function $\phi_t$ can be, e.g., a linear, quadratic, \blue{or piecewise linear} function.}

%

\subsection{Non-convex scheduling problem}\label{scheduling:aggregator}
The aggregator aims to identify an admissible schedule for each job while minimizing the overall cost of serving the resulting aggregate load. 
This scheduling problem can be cast as a mixed-integer optimization as follows. 

Let $s_j(t)$, $j\in \Jcal$ and $t\in \Tcal$, be the indicator of job starting times: 
\begin{equation}
s_j(t) = \begin{cases}
	1, & \quad \mbox{if job $j$ starts in period $t$},\\
	0, & \quad \mbox{otherwise.}
\end{cases}\nonumber
\end{equation}
We have
\begin{align}
s_j(t) \in \{0,\,1\}, &\quad j\in \Jcal,\,\, t\in \Tcal,\nonumber\\
 \sum_{t\in \Tcal}s_j(t) = 1, &\quad j \in \Jcal.\nonumber
 \end{align}

For a schedule of job $j$, embedded by \blue{$\mathbf s_{j}$}, to be admissible, we \blue{also }require
\begin{align}
s_j(t) =0, \quad j \in \Jcal, \,\,  t \not\in \Tcals_j.\nonumber 	
\end{align}

 Define matrix \blue{$\mathbf P^{(j)}\in \reals^{T\times T}$}, $j \in \Jcal$, to be a dictionary of power consumption profiles for job $j$, when its starting time varies. 
 In other words, the $t$-th column of \blue{$\mathbf P^{(j)}$} is the power consumption profile of job $j$ with starting time $t$, i.e., for $t,\, t'\in \Tcal$,
 \begin{equation}\label{eq:def:Pj}
\orange{P^{(j)}_{t', t}} = \begin{cases}
	p_j, & \quad \mbox{if } t\le t' \le t+d_j-1,\\
	0, & \quad \mbox{otherwise.}
\end{cases}
 \end{equation}
The power consumption profile of job $j$ with starting time embedded in \blue{$\mathbf s_{j}$} is then 
\begin{equation}
\bm \ell_j = \mathbf P^{(j)} \mathbf s_j \in \reals^T, \quad j \in\Jcal.\nonumber
 \end{equation}
 \rev{Indeed, if $s_j(t)=1$ and $s_j(t')=0$ for all $t'\neq t$, then \blue{$\bm \ell_j$} calculated as above coincides with the $t$-th column of $\mathbf P^{(j)}$.}
 
 Thus, the scheduling optimization takes the following form
\begin{subequations}\label{opt:s}
 	\begin{align}
 		\min_{\mathbf s \in \reals^{J\times T},\,\, \mathbf  L \in \reals^T} \quad & \Phi(\mathbf L)\\
 		\mbox{s.t.} \qquad\quad & \mathbf L = \sum_{j\in \Jcal}  \mathbf P^{(j)} \mathbf s_j,\label{opt:s:a} \\
 		& s_j(t) \in \{0,\,1\}, \quad j\in \Jcal,\,\, t\in \Tcal,\label{opt:s:b}\\
 		 &  s_j(t) =0, \qquad\quad j \in \Jcal,\,\, t \not\in \Tcals_j, \label{opt:s:c}\\
 &\sum_{t\in \Tcal}s_j(t) = 1, \label{opt:s:d} \quad\, j \in \Jcal.
 	\end{align}
 \end{subequations}
Problem \eqref{opt:s} is challenging due to the potential problem size (i.e., large $J$) and the binary constraint \eqref{opt:s:b}.
In Section~\ref{sec:scheduling}, we propose an efficient algorithm for solving the problem with provable performance guarantees.
\def\Ecal{\mathcal E}
\section{\blue{Relaxation-Adjustment-Rounding Algorithm}}\label{sec:scheduling}
 \blue{In this section,} we propose a scalable procedure for obtaining a near-optimal solution to the mixed-integer convex program \eqref{opt:s} in polynomial time, and we will refer to this procedure throughout the paper as \blue{the} \blue{\emph{relaxation-adjustment-rounding}} \blue{(RAR)} algorithm.  
The proposed procedure consists of three major steps, which are detailed as follows. 

\subsubsection{Relaxation} 
We start by solving a convex relaxation to the mixed-integer convex program \eqref{opt:s}, where the binary constraint~\eqref{opt:s:b} is replaced by 
\begin{equation}
	s_j(t) \in [0,1], \quad j\in \Jcal,\,\, t\in \Tcal.\label{opt:s:bc}
\end{equation}
\purple{The resulting optimization problem is a convex program that can be solved in polynomial time.}
 We denote the obtained schedule of this step by $\mathbf s^\mathrm{R}$ and its associated aggregate load profile by $\mathbf  L^{\mathrm{R}}$. 

\subsubsection{Adjustment} 
The solution obtained in the previous step may contain many fractional entries. 
The goal of this step is to reduce the number of fractional entries in  $\mathbf s^\mathrm{R}$ to a  sufficiently small number (see Lemma \ref{lemm:4} for the exact upper bound) without changing the aggregate load. 

To this end, we consider the following undirected \emph{multigraph} $G_d(\blue{ \mathbf s}) = (\Tcal, \Ecal_{d})$ for each duration $d =1, \dots, d^{\max}$, \blue{defined for any $\mathbf s$}. \blue{As we shall see, we will construct the graph first for $\mathbf s =\mathbf s^\mathrm{R}$. }
The node set of the multigraph is the set of time periods.  
We create the edge set of the multigraph by adding one edge for each $(t',\,\,t,\,\,j)\in \Tcal \times \Tcal \times \Jcal_d$ such that  
\begin{equation}\label{eq:edge:multi}
t' = \min \{\blue{\tilde t}:s_j(\blue{\tilde t}) \not\in \{0,1\}\}, \ s_j(t) \not\in \{0,1\}, \ t' \neq t, 
\end{equation}
where $\Jcal_d$ is the set of jobs with duration $d$. 
In other words, for any job $j$, we find the first time $t'$ such that $s_j(t')$ is fractional and then add an edge $(t',\,t,\, j)$ to the edge set $\Ecal_{d}$ for each $t\neq t'$ with a fractional schedule $s_j(t)$. 
It is therefore convenient to refer to the edges in the multigraph by $(t',\,\,t,\,\,j)$, whose first two indices are the two nodes connected by the edge, and the last index can be viewed as a label to distinguish multiple edges connecting the pair of nodes. 
For any edge $e\in \Ecal_{d}$, we denote the corresponding \orange{triple} by $(t'_e,\,\, t_e,\,\, j_e)$. 

Given the definition of $G_d(\blue{\mathbf s})$, we can summarize the adjustment step in Algorithm~\ref{alg:adj}. 
For each $d=1, \dots, d^{\max}$, the algorithm starts by constructing $G_d(\mathbf s^\mathrm{R})$. 
It then tries to find a cycle in the multigraph, identify an adjustment that eliminates at least one fractional entry in $\mathbf s^\mathrm{R}$ (and therefore also one edge in the multigraph), and update the schedule and graph accordingly until no further cycles can be found. We denote the result of this adjustment process by $\mathbf s^\mathrm{A}$ and its associated aggregate load profile by $\mathbf{L}^{\mathrm{A}}$. As we shall see later in Section~\ref{sec:analysis}, this adjustment process is lossless, i.e., $\Phi(\mathbf{L}^\mathrm{A})=\Phi(\mathbf{L}^\mathrm{R})$. 
Therefore, \emph{$\mathbf{s}^\mathrm{A}$ is another optimal point of the relaxed problem, with no more fractional entries than $\mathbf{s}^\mathrm{R}$}.

\begin{algorithm}[tp]
    \caption{Lossless adjustment, rectangular loads}\label{alg:adj}
    \label{alg::1}
    \LinesNumbered
    \KwIn{$\mathbf s\gets  \mathbf s^\mathrm{R}$\purple{.}}
  \For{$d=1,\,\,\dots,\,\,d^{\max}$}{
  \While{$G_d(\blue{\mathbf s})$ is cyclic}{
  Find a cycle $\mathcal{C}$ in $G_d(\blue{\mathbf s})$
  (Treat edges in $\mathcal C$ as directed such that the directed edges still form a cycle)\purple{.}\\
  Find $\Delta^\star$ by solving:
		\begin{subequations}\label{cons:alg}
		\begin{align}
		\max_{\Delta}\quad & \Delta \\
		\mbox{s.t.} \quad 
		& s_{j_e}(t'_e)-\frac{\Delta}{p_{j_e}}\ge 0, \quad   e\in \mathcal{C}\purple{,}  \label{cons:1}\\
		& s_{j_e}(t_e)+\frac{\Delta}{p_{j_e}}\le 1, \quad   e\in \mathcal{C}\purple{.} \label{cons:2}
		\end{align}
		\end{subequations}
		\\
		Update schedule $\mathbf s$:	
		\begin{subequations}\label{cons:upd}
		\begin{align} 
		& s_{j_e}(t'_e)\gets s_{j_e}(t'_e)-\frac{\Delta^\star}{p_{j_e}},\quad   e\in \mathcal{C} \purple{,}  \label{cons:3}\\
		& s_{j_e}(t_e)\gets s_{j_e}(t_e)+\frac{\Delta^\star}{p_{j_e}},\quad   e\in \mathcal{C}\purple{.} \label{cons:4}
		\end{align}	
		\end{subequations}
		\\
		Update $G_d(\blue{\mathbf s})$ with the new schedule $\mathbf s$.
	}
  }
  \KwOut{$\mathbf s^\mathrm{A}\gets \mathbf s$\purple{.}}
\end{algorithm}
	\begin{figure*}[h]
		\centering
		\includegraphics[width=1\textwidth]{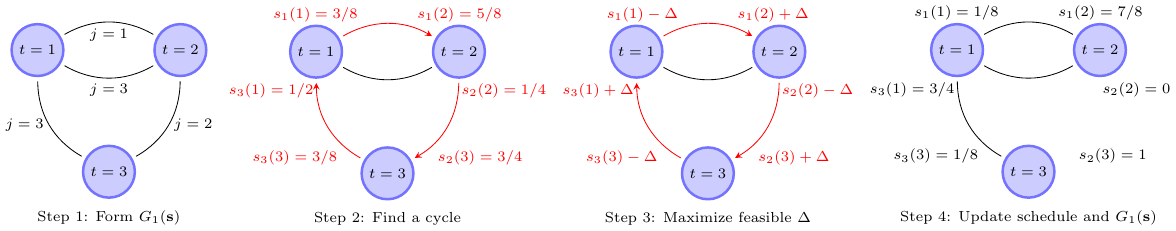}
		\caption{Forming and updating the multigraph for the adjustment step: A toy example}\label{alg:example}
	\end{figure*}
\begin{example}[Adjustment for rectangular load shape]
	\brown{We present a simple example  to show how to form the graph and perform an iteration of the while loop in Algorithm~\ref{alg:adj}. Consider a setting where $T=3$ and there are different jobs with various durations but four of them have duration $d=1$. \purple{Let $p_{j}=1$ for all of these four jobs.} After  the relaxation step, suppose that the obtained schedule for these jobs is as follows:
		\begin{align}
			\begin{bmatrix}
				{s}_{1}^{\mathrm{R}}(1) & {s}_{1}^{\mathrm{R}}(2)&{s}_{1}^{\mathrm{R}}(3)\\[2pt]
				{s}_{2}^{\mathrm{R}}(1) & s_{2}^{\mathrm{R}}(2) & {s}_{2}^{\mathrm{R}}(3)\\[2pt]
				{s}_{3}^{\mathrm{R}}(1) & {s}_{3}^{\mathrm{R}}(2) & {s}_{3}^{\mathrm{R}}(3)\\[2pt]
				{s}_{4}^{\mathrm{R}}(1) & 
				{s}_{4}^{\mathrm{R}}(2) &{s}_{4}^{\mathrm{R}}(3)
			\end{bmatrix} &= \begin{bmatrix}
				\frac{3}{8} & \frac{5}{8} & 0\\[2pt]
				0 & \frac{1}{4} &  \frac{3}{4}\\[2pt]
				\frac{1}{2} & \frac{1}{8} & \frac{3}{8}\\[2pt]
				0 &1&0
			\end{bmatrix}.\label{alg1:example}
		\end{align}
		Fig.~\ref{alg:example} illustrates the steps.} \purple{In \gray{Step 1}, we add edges associated with pairs of fractional entries in different time slots to form $G_{1}\left(\mathbf{s}\right)$. \gray{For example, given $s_{1}^{\mathrm{R}}(1)$ and $s_{1}^{\mathrm{R}}(2)$ are fractional, we add an edge between node $t=1$ and $t=2$, labeled by $j=1$, i.e., edge $(1,2,1)$ in our \orange{triple} notation.} We then find a cycle in Step 2 and proceed to Step 3 by performing~\eqref{cons:alg} to find $\Delta^\star$. This parameter is then used in Step 4 to update the variables associated with the endpoints of the edges in the selected cycle according to~\eqref{cons:upd}. As a result, fractional entries $s_{2}(2)$ and $s_{2}(3)$  are converted into \gray{integers} and then $G_1(\mathbf s)$ is updated accordingly. One can verify that the aggregate load of these four jobs is not changed after doing this adjustment.}  
		\end{example}
\subsubsection{Rounding} 
The result of the adjustment step may still contain fractional entries. 
We round the elements in $\mathbf{s}^\mathrm{A}$ by scheduling job $j$ to start at time $t$ with probability $s^\mathrm{A}_{j}(t)$ for each $j\in \Jcal$ and $t\in \Tcal$ such that $s_j^\mathrm{A}(t)>0$. 
We denote the resulting integral schedule by $\mathbf{s}^\mathrm{I}$, and its associated aggregate load profile by $\mathbf{L}^\mathrm{I}$. 
\rev{Note that even though this step introduces randomness to our algorithm, our analysis is performed based on bounding the number of fractional entries in $\mathbf{s}^\mathrm{A}$. Therefore, our results hold for any realization of this randomized rounding step and thus deterministically.}

\purple{\section{Performance of the RAR algorithm}\label{sec:IV}}
\purple{In this section, we start by analyzing the RAR algorithm and then briefly discuss the necessity of the adjustment step. Detailed proofs can be found in the appendices.}
\subsection{Analysis of the RAR algorithm}\label{sec:analysis}
Our main result regarding the RAR algorithm is that it is computationally efficient and it provides a near-optimal (in fact asymptotically optimal in a per-job cost sense) solution. We start by analyzing the computational complexity of the proposed algorithm. Since the convex relaxation step and the rounding step finish in polynomial time, we focus on the adjustment step. Proofs omitted in the main text can be found in the appendices. 


\begin{lemma}[Complexity of the adjustment step]\label{lemm:1}
The while loop in Algorithm~\ref{alg:adj} terminates in  $O(J T)$ iterations across all $d$ values, with the complexity of each iteration being $O(T)$.  
\end{lemma}
 
\purple{Given Lemma~\ref{lemm:1}, the overall complexity of the adjustment step is $O(JT^2)$, which is linear \rev{in} the number of jobs.} 
When the availability windows of all jobs are considerably smaller than $T$, it is not hard to show that our complexity bound can be tighten to $O(J T\tau^{\max})$, where $\tau^{\max} = \max_{j\in \Jcal} \tau_j$. 

To analyze the performance of the proposed algorithm, we denote an optimal point of the  mixed-integer program by $\mathbf{s}^\star$ and its associated aggregate load profile by $\mathbf{L}^\star$. 
Then, we have:
\begin{theorem}[Sub-optimality bound]\label{theo:1}
	The schedule $\mathbf{s}^\mathrm{I}$ is feasible for \blue{the mixed integer program}~\eqref{opt:s}, with sub-optimality bounded as follow\footnote{\gray{We can tighten the bound in~\eqref{bounds} by a factor of 2 if $\phi_t(\cdot)$ is non-decreasing for all $t \in \mathcal{T}$.}}:
	\begin{equation}
		\Phi(\mathbf{L}^\mathrm{I}) - \Phi(\mathbf{L}^\star) \le \purple{2}d^{\max}TK \max_{j\in \mathcal{J}}p_jd_j.\label{bounds}
	\end{equation}
\end{theorem}

In other words, the schedule produced by our RAR algorithm is near optimal, with a sub-optimality bound that depends on the job characteristics, the length of the decision horizon, and the Lipschitz coefficient of the cost functions. 
Notably, the sub-optimality bound is independent of the number of jobs $J$. 
As a consequence, the per-job sub-optimality, computed in~\eqref{eq:perjob} below, approaches zero when the size of the job set increases:
\textcolor{black}{
	\begin{corollary}[Asymptotic optimality for per-job cost]\label{cor:Th:1}
		For any $J$ jobs with a bounded $\max_{j\in \Jcal} p_j d_j$,
		\begin{equation}\label{eq:perjob}
			\frac{\Phi(\mathbf{L}^\mathrm{I}) - \Phi(\mathbf{L}^\star)}{J} \to 0, \quad \mbox{ as } J\to\infty. 
		\end{equation}  	
\end{corollary}}

\blue{While the detailed proofs are given in the appendices, we outline the key steps in establishing our performance theorem next, which will  elaborate on the design principle of the adjustment step. }
\rev{The proof of Theorem~\ref{theo:1} relies on the following three key lemmas.} 
\begin{lemma}[Feasibility of the \purple{adjusted} solution]\label{lemm:2}
	The schedule $\mathbf{s}^\mathrm{A}$ satisfies constraints~\eqref{opt:s:c} and~\eqref{opt:s:d}, and therefore is feasible for the convex relaxation of~\eqref{opt:s}. 
\end{lemma}
	\begin{lemma}[Lossless adjustment]\label{lemm:3}
	The adjustment step is lossless, i.e., 
	\begin{equation}
	\Phi(\mathbf{L}^\mathrm{A}) = \Phi(\mathbf{L}^\mathrm{R})  \le \Phi(\mathbf  L^\star).\nonumber
	\end{equation}
	\end{lemma}

By Lemma~\ref{lemm:2} and Lemma~\ref{lemm:3}, we observe that the adjustment step takes an optimal point of the convex relaxation $\mathbf{s}^\mathrm{R}$, and produces another feasible schedule $\mathbf{s}^\mathrm{A}$ with the same cost. 
	It follows that $\mathbf{s}^\mathrm{A}$ is another optimal point of the convex relaxation. 
	Furthermore, we can characterize the number of fractional entries in $\mathbf{s}^\mathrm{A}$ as follows.

	\begin{lemma}[Bound on number of fractional entries]\label{lemm:4}
		The number of fractional entries in $\mathbf{s}^{\mathrm{A}}$ is upper bounded by $2d^{\max}T$.
	\end{lemma}
	
	This allows us to bound the loss of optimality introduced by the rounding step and eventually characterize the overall sub-optimality of the algorithm. 

\subsection{Necessity of the adjustment step}
\blue{An intuitive folklore is that the problem of non-preemptive scheduling may be well approximated by a continuous optimization problem when the number of jobs is large. However, we are unaware of any concrete prior algorithm or analysis demonstrating this can be done in a rigorously provable manner. As the relaxation and rounding steps in the RAR algorithm are also intuitive (and may have been attempted in the past), one interesting question is \emph{whether the adjustment step is necessary in achieving the performance guarantee in Theorem~\ref{theo:1}}? Our next result gives an affirmative answer to this question by showing that the intuitive relaxation-rounding algorithm (that skips the adjustment step) can lead to arbitrarily bad performance.  }

    \begin{lemma}[Necessity of adjustment step] \label{adjustment:necess}
    	For any natural number $N$, there exist a problem instance and a non-zero probability\footnote{The statement is probabilistic as our rounding scheme is so. It is not hard to show that the same observation also holds for the deterministic rounding scheme that always schedules a  job with fractional $\mathbf s_j^\mathrm{R}$ to start at the time with the maximum $s^\mathrm{R}_j(t)$. } such that the relaxation-rounding algorithm that directly rounds $\mathbf s^\mathrm{R}$ into $\mathbf s^\mathrm{I}$ (with the corresponding aggregate load denoted by $\mathbf L^\mathrm{I}$) satisfies:
    	 \begin{equation}\label{adjustment:necess:eq}
     \frac{1}{J} \left[\Phi\left(\mathbf{L}^{\mathrm{I}}\right) - \Phi\left(\mathbf{L}^\star\right) \right] \ge \frac{1}{2} NT. 
    \end{equation}
    \end{lemma}

    \blue{This result states that the  performance of the RAR algorithm without the adjustment step can be arbitrarily bad: the per-job sub-optimality gap, which approaches $0$ for the RAR algorithm, can be $\Omega(NT)$ for an arbitrary $N$ in this case. 
    The stark contrast
     highlights the necessity of the adjustment step. It also speaks to the importance of algorithm design -- without developing our adjustment step, the folklore as it stands lacks a proper theoretical ground.  }

\section{Extensions}\label{sec:ext}
\subsection{Incorporating uncertainties}\label{stochastic}
\blue{It is critical to model and incorporate uncertainties for real world implementation of the scheduling scheme that we develop in this paper. There may be uncertainties associated with the job parameters, distributed renewables that may be used by the aggregator to serve the loads locally, and/or electricity market prices if the aggregate load is served by purchasing electricity from the grid. While the focus of this paper is on addressing the challenges associated with the non-convexity of the scheduling problem induced by binary variables, we here extend our algorithm and results to a class of stochastic instances of the problem of the following form:}
\begin{subequations}\label{opt:s:u}
 \begin{align}
 		\min_{\mathbf s \in \reals^{J\times T},\,\, \mathbf  L \in \reals^T} \quad & \Phi(\mathbf L):=\mathbb E_{\bm \theta} \left[\sum_{t\in \Tcal} \mathbb  \phi_{t}\left(L(t);\bm{\theta}(t)\right)\right] \label{opt:s:u:obj}\\
 	\mbox{s.t.} \qquad\quad & \eqref{opt:s:a},\eqref{opt:s:b},\eqref{opt:s:c},\eqref{opt:s:d}\orange{,}
  	\end{align} 
  \end{subequations}%
where $\bm \theta(t)$ is a random vector including all the uncertain parameters under consideration, and we have overloaded the notation $\phi_t$ to model the stage-wise cost function that may \orange{depend} on the uncertain parameters. Notably, \eqref{opt:s:u} only addresses uncertainties that impact the scheduling cost, uncertainties such as job parameters that affect the constraints are not treated and are left for future work. 

The  stochastic program \eqref{opt:s:u} can be readily solved by our proposed RAR algorithm, provided that the expectation in the cost function can be efficiently evaluated. 
Indeed, given the distribution of $\bm \theta$, one can evaluate the expectation using, e.g., Monte Carlo simulation, carry out the RAR algorithm as it is, and obtain the same performance guarantee in Theorem~\ref{theo:1} (adjusted by errors due to Monte Carlo). 

There are also occasions when the expectation in \eqref{opt:s:u} may be difficult to evaluate. This can be the case when  the (potentially high-dimensional) joint distribution  of the uncertain parameters $\bm \theta(t)$ is not available. As such, the inclusion of uncertainties in the scheduling cost introduces new challenges as we may not be able to efficiently and accurately evaluate the cost function.   

We proceed to demonstrate that this challenge can be circumvented by leveraging the structure of the  problem at hand, when the cost functions \gray{satisfy the following two assumptions}:

\begin{assumption}[Time-invariant costs]\label{assum:1} There exists a \gray{time-invariant} strongly convex and $K_{\psi}$-Lipschitz continuous \gray{differentiable} function $\psi(\cdot): \reals\mapsto \reals$ \purple{with $K_{\psi}<\infty$} such that
	 \begin{equation}
			\mathbb E_{\bm \theta} \left[\sum_{t\in \Tcal} \mathbb  \phi_{t}\left(L(t);\bm{\theta}(t)\right)\right]=\sum_{t\in \Tcal}\psi(L(t)).\label{cost:uncertain}
		\end{equation}
\end{assumption}%

%
%
%
%

The assumption only states the existence of $\psi$; in other words, we do not assume the function $\psi$ is known or can be efficiently evaluated. This assumption holds under some settings. For example, when the cost functions $\phi_t$, $t\in \mathcal T$, do not vary with $t$ and are strongly convex, and $\{\bm \theta(t): t\in \mathcal T\}$ is a stationary stochastic process, it is easy to verify that Assumption~\ref{assum:1} holds. The next example gives another setting where \purple{Assumption~\ref{assum:1} holds without $\phi_t$ being strongly convex}:
{\begin{example}[Load serving with renewables]
\blue{Consider the setting where the aggregator serves a portion of the aggregate load with distributed renewable generation. In this case, we may use $\bm \theta(t)$ to denote the renewable outputs from different generation units.
	  The excess load may be served by purchasing electricity from the utility \purple{backstop supply}, so the cost function takes the form of $\phi_t(L(t); \bm \theta(t))= \alpha (L(t)- \mathbf 1^\top \bm \theta(t))_+$, where $\alpha$ is the utility electricity price. If the aggregator decides to model $\{\bm \theta(t): t\in \mathcal T\}$ as a stationary process (e.g., when distributed wind generation is used and there lack temporally granular data to fit a time-varying model), \purple{such that $\bm \theta(t)$ has a density function,} Assumption~\ref{assum:1} holds\cite{bertsekas1973stochastic}.
}
\end{example}


\begin{assumption}[Generalized monotonicity]\label{assum:2}
	The following property holds for the derivative of function $\psi$. For any $d=1, \dots, d^{\max}$,  $\mathbf x \in \mathbb R^d$, and $\mathbf y \in \mathbb R^d$, we have 
	\begin{equation}\label{eq:gen:monotone}
		\mathbf 1^\top \psi'(\mathbf x) < \mathbf 1^\top \psi'(\mathbf y) \mbox{ if and only if } \mathbf 1^\top \mathbf x < \mathbf 1^\top \mathbf y,
	\end{equation}
	where we slightly abuse the notation and use $\psi'(\mathbf x)$ to denote $\{\psi'(x_i)\}_{i=1}^d$. 
\end{assumption}

This assumption states that $\psi'$ satisfies a notion of monotonicity that connects the sum of  inputs and the sum of their function values. When $d^{\max}=1$, this assumption is implied by the strong convexity of $\psi$ as \eqref{eq:gen:monotone} is equivalent to the strict monotonicity of $\psi'$. In contrast, when $d^{\max}>1$, $\psi'$ being strictly increasing is  a necessary but not sufficient condition for \eqref{eq:gen:monotone}. A sufficient condition, which may hold for certain practical settings, is that $\psi'$ is an increasing affine function (or $\psi(z) = a z^2 + bz$ for some $a>0$ and $b$). 
%

\blue{Under Assumptions~\ref{assum:1} and~\gray{\ref{assum:2}}, it turns out that we can solve~\eqref{opt:s:u} by applying a slightly modified version of the proposed RAR algorithm \emph{without evaluating the expectation or function $\psi$.} In this modified procedure, we simply replace the cost function $\Phi(\mathbf L)$ with $\purple{\widetilde{\Phi}}(\mathbf L):= \sum_{t\in \mathcal T} L(t)^2$ and apply the RAR algorithm with the new cost function. Our next theorem states that the schedule obtained via this procedure enjoys a similar performance guarantee as in Theorem~\ref{theo:1}.  }

\blue{\begin{theorem}[Sub-optimality bound, modified RAR]\label{uncertain:result}
Under Assumptions~\ref{assum:1} and~\ref{assum:2},  the obtained scheduling of the modified RAR algorithm, $\mathbf{L}^{\mathrm{I}}$, is feasible for~\eqref{opt:s:u} with sub-optimality bounded as follows:
\begin{align}
	\Phi(\mathbf{L}^\mathrm{I}) - \Phi(\mathbf{L}^\star) 
	\leq \purple{2}d^{\max}TK_{\psi} \max_{j\in \mathcal{J}}p_jd_j.\label{un:bounds}
\end{align}
\end{theorem}}
\blue{Theorem~\ref{uncertain:result} is established by leveraging a \emph{cost equivalence property} of the convex relaxation of~\eqref{opt:s:u} when time-invariant costs satisfying Assumptions~\ref{assum:1} and~\ref{assum:2} are considered (\purple{refer to} Lemma~\ref{th:cost} in Appendix~\ref{appB} for details): The set of optimal points for the convex relaxation of~\eqref{opt:s:u} and that of the same optimization with the objective function replaced by $\purple{\tilde{\Phi}}$ is identical under Assumptions~\ref{assum:1} and~\ref{assum:2}.  Equipped with this property, we can obtain an optimal point of the convex relaxation of \eqref{opt:s:u} without evaluating its actual cost function, and then apply the adjustment and rounding steps to obtain a feasible and near-optimal solution for \eqref{opt:s:u}. 
}
 

\subsection{\blue{Decentralization}}\label{sec:pricing} 
\blue{So far we have focused on the centralized setting where the aggregator or load-serving entity \emph{directly} schedules the jobs. In practice, unless the users surrender the control of their loads to the aggregator (e.g., by participating in certain demand response programs), a decentralized setting where the aggregator \emph{indirectly} control the loads via incentives/prices may be more appropriate. }

\blue{In this section, we provide a limited account of this decentralized setting by studying properties of a pricing policy for indirect load control that is natural for our problem. In particular, we consider the \emph{marginal pricing} policy defined for the convex relaxation of~\eqref{opt:s}. Let 
} $\bm \lambda\in \reals^T$ be the dual solution associated with the constraint~\eqref{opt:s:a} in the relaxed version of problem~\eqref{opt:s} where the binary constraint \eqref{opt:s:b} is replaced with the convex constraint \eqref{opt:s:bc}. 
\blue{A user facing the marginal pricing policy who chooses \purple{any} admissible schedule $\mathbf s_j'$ with induced load $\bm \ell_j'$ needs to pay 
\begin{equation}
	\pi_j^\mathrm{MP} (\mathbf{s}_j';\,\, \bm\lambda)= \sum_{t\in \Tcal} \lambda(t) \ell'_j(t) = \bm \lambda^\top \mathbf P^{(j)}\mathbf s_j',\quad j\in \Jcal. \label{eq:mpp}
	\end{equation}
	Note that the pricing policy~\eqref{eq:mpp} applies regardless whether the user decided schedule $\mathbf s_j'$ coincides with the schedule $\mathbf s_j^\mathrm{I}$ computed by the RAR algorithm.
}

\blue{The primary research question that is of concern in this decentralized setting is \emph{whether the marginal pricing policy~\eqref{eq:mpp} will incentivize a ``good'' schedule as evaluated by the total cost for serving the aggregate load $\Phi$.} We do not expect that we may achieve an \emph{optimal} schedule  defined by the solution of mixed integer program~\eqref{opt:s}, as fundamentally there is a gap between the mixed integer program and its convex relaxation, \purple{based on which the marginal prices are defined}. However, somewhat surprisingly, we can establish that the marginal pricing rule can indeed incentivize the schedules produced by the RAR algorithm, which by Theorem~\ref{theo:1} is a near-optimal solution of~\eqref{opt:s}. 
This offers another evidence that our RAR algorithm and the theoretical results around it can serve as a bridge formally connecting the original mixed integer program and its convex relaxation.} \purple{To simplify the exposition, we introduce the following additional assumption for results in this subsection.}

\blue{\begin{assumption}[Differentiability and constraint qualification]
	The cost function $\Phi$ is differentiable in its domain. Furthermore, the Linear Independent Constraint Qualification (LICQ) holds for the convex relaxation of~\eqref{opt:s}. 
\end{assumption}}

\blue{Given the differentiability of the cost function, by the optimality condition of the convex relaxation, it is easy to show that 
\begin{equation}
	\lambda(t) = \frac{\partial \Phi(\mathbf{L}^\mathrm{R})}{\partial L^\mathrm{R}(t)} = \frac{\mathrm{d} \phi_t(L^\mathrm{R}(t))}{\mathrm{d}L^\mathrm{R}(t)}, \quad t\in \Tcal,\nonumber
\end{equation}
hence $\lambda(t)$ is indeed the marginal price in time period $t$. Meanwhile, the LICQ is a mild condition that is commonly assumed to rule out degeneracy and will ensure the uniqueness of the dual solutions \cite{wachsmuth2013licq}  (and hence the marginal prices) for the convex relaxation of~\eqref{opt:s}.
}


\blue{\purple{The first} step in bridging the gap between a payment rule defined based on the convex relaxation and the sub-optimal integer solution to~\eqref{opt:s} produced by the RAR algorithm is the following lemma.
}

\begin{lemma}[\blue{Payment equivalence}]\label{lem:price}
	For every job $j\in \mathcal{J}$,  \blue{the payments calculated from the marginal pricing rule~\eqref{eq:mpp} with an optimal relaxed schedule $\mathbf s^\mathrm{R}_j$ and with the schedule produced by the RAR algorithm $\mathbf s^\mathrm{I}_j$ are identical:}
	\begin{equation}
	\pi_j^\mathrm{MP}(\mathbf s^\mathrm{R}_j;\,\, \bm \lambda)=\pi_j^\mathrm{MP}(\mathbf s^\mathrm{I}_j;\,\, \bm \lambda), \blue{\quad j \in \mathcal J}.\nonumber
	\end{equation}
\end{lemma}

\blue{Lemma~\ref{lem:price} is established leveraging the structural properties of the primal and dual solutions of the convex relaxation, and how the RAR algorithm obtains an admissible schedule from the fractional solution $\mathbf s^\mathrm{R}$. Equipped with this key lemma, we can establish our main result on the decentralization of the schedules produced by the RAR algorithm using marginal prices:}

\begin{theorem}[Self-scheduling]
\label{SS}
	For each job $j\in \Jcal$, and any admissible $\mathbf s_j'$, there is   
	\begin{equation}
		\pi_j^\mathrm{MP} (\mathbf s_j^{\mathrm{I}};\,\, \bm \lambda) \le \pi_j^\mathrm{MP} (\mathbf s_j';\,\, \bm \lambda), \quad j \in \Jcal. \label{no:incntive}
	\end{equation} 
\end{theorem}

Theorem~\ref{SS} states that if job owners are charged \blue{according to} the marginal pricing rule~\eqref{eq:mpp}, for each job $j\in \Jcal$, there is no incentive to unilaterally deviate from the schedule $\mathbf s_{j}^{\mathrm{I}}$ computed by the aggregator. 
\blue{In other words, the marginal price $ \bm \lambda$ incentivizes decentralized implementation of the RAR solution $\mathbf s^{\mathrm{I}}$. }

\begin{remark}[\blue{Weakly dominant strategy}]
\blue{The weak inequality in~\eqref{no:incntive} suggests that given the marginal price $\bm \lambda$, implementing $\mathbf s^\mathrm{R}_j$ is a \emph{weakly dominant strategy}. In other words, there may be other schedules  that leads to an identical cost for job $j$. Were these alternative schedules picked by the jobs, the aggregate load may not be \gray{$\mathbf{L}^\mathrm{I}$} and our performance guarantees for the RAR algorithm may not apply. Therefore, if the aggregator utilizes the marginal prices for a decentralized implementation, the recommended schedule $\mathbf s^\mathrm{I}_j$ should be communicated to job $j$ in addition to the prices $\bm \lambda$. }
\end{remark}

\begin{remark}[Revenue adequacy and flexibility-revealing properties]\label{rk:pricing:etc}
	In addition to self-scheduling, we establish in \cite{chen2022scheduling} that the marginal pricing policy is approximately revenue adequate and can incentivize the truthful reporting of the flexibility window of jobs \emph{when there is only single-slot jobs}, i.e., when $d^{\max} = 1$. When $d^{\max}>1$, this flexibility-revealing property may not hold. 
	Thus\orange{,} while the marginal pricing policy can incentivize the RAR solution, it does not satisfy some of the desirable properties for a good pricing policy. 
\end{remark}

\subsection{Realistic load shapes}\label{sec:realload}
In practice, the load shapes are rarely rectangular. In this section, we generalize our algorithm and results to the case with realistic load shapes \purple{where the power consumption of jobs vary across time slots.}

To model non-preemptive jobs with realistic load shapes, we overload our notation and refer to the time-varying power requirement of the load $j$ by $\mathbf{p}^{(j)}\in \reals^{d^{\max}}$, where
$p^{(j)}_i \ge 0$ for $i= 1,\dots, d_j$ and $p^{(j)}_i=0$ for $i = d_j+1, \dots,  d^{\max}$.   
Moreover, the dictionary of power consumption profiles for job $j$ \brown{can be denoted by $\mathbf P^{(j)}\in \reals^{T\times T}$ in which its $t$-th column denotes the time-varying power consumption of the job if it starts at time $t$. In other words, for every $t, t' \in\mathcal T$,}
\orange{\begin{equation}
P^{(j)}_{t', t} = \begin{cases}
	p^{(j)}_{ t'-t+1}, & \quad \mbox{if } t\le t' \le t+d_j-1,\\
	0, & \quad \mbox{otherwise.}\nonumber
\end{cases}
\end{equation}}

\subsubsection{Scheduling algorithm}
The scheduling algorithm follows the same three-step procedure as outlined in Section~\ref{sec:scheduling}, with the adjustment step modified to incorporate non-rectangular loads.
In particular, Algorithm~\ref{alg:adj} cannot be directly applied for non-rectangular loads. \purple{This is because} when two jobs have different shapes, we cannot use the increment of the schedule of a job, \purple{i.e., certain fractional entries,} to balance the change in total power consumption caused by the reduction of the schedule of the other job. As a result, \gray{Lemma~\ref{lemm:3}} does not hold if we directly apply Algorithm~\ref{alg:adj} to jobs with general load shapes.
\purple{To address this challenge, we propose a new adjustment procedure summarized in Algorithm~\ref{alg:adjArb}.} 


Let $D=d^{\max}+1$. The algorithm starts by finding \purple{a time pair $t, \tilde t$ such that there exists $D$ jobs to form $\mathcal{J'}:=\{j_1,\dots, j_{D}\}\subset \mathcal{J}$ such that for each $j_k\in \mathcal{J'}$, }
%
\begin{equation}\label{pair}
s_{j_k}(t)\not= \{0, 1\},\quad  s_{j_k}(\tilde t\,)\not= \{0, 1\},\quad t\neq \tilde t.
\end{equation}
Consider the matrix $\mathbf P_{\mathcal{J'}}=[ \mathbf p^{(j_1)},\dots, \mathbf p^{(j_D)}] \in \mathbb R^{(D-1)\times D}$, whose rank is less than $D$. We can then find a non-zero element in the null space of $\mathbf P_{\mathcal{J'}}$, i.e., a non-zero vector $\bm \xi \in \mathbb R^{D}$ such that
\begin{equation}\label{eq:xi}
	\mathbf P_{\mathcal{J'}} \bm \xi = \mathbf 0. 
\end{equation}
Using $\bm \xi$, similar to Algorithm 1, the algorithm identifies an adjustment that eliminates at least one fractional entry in $\mathbf s^\mathrm{R}$, and updates the schedule until no further jobs and time pair can be found. 
\begin{example}[Adjustment for realistic load shape]
Suppose $T=4$ and there are only 4 jobs with $d_{1}=1$ and \purple{$d_{2}=d_{3}=d_{4}=2$}. Thus, $d^{\mathrm{max}}=2$. The time-varying power requirement of the jobs are $\mathbf{p}^{(1)}= \begin{bmatrix}
		1, 0 
	\end{bmatrix}^\top$, $	\mathbf{p}^{(2)}= \begin{bmatrix}
	1, 2
	 
	\end{bmatrix}^\top$,  $	\mathbf{p}^{(3)}= \begin{bmatrix}
	1, 3
	 
	\end{bmatrix}^\top$,  $	\mathbf{p}^{(4)}= \begin{bmatrix}
	2, 3
	 
	\end{bmatrix}^\top$. Suppose after the relaxation step, the following schedule is obtained for these jobs:
 \begin{align}
	\begin{bmatrix}
		{s}_{1}^{\mathrm{R}}(1) & {s}_{1}^{\mathrm{R}}(2)&{s}_{1}^{\mathrm{R}}(3)&{s}_{1}^{\mathrm{R}}(4)\\[2pt]
		{s}_{2}^{\mathrm{R}}(1) & s_{2}^{\mathrm{R}}(2) & {s}_{2}^{\mathrm{R}}(3)&{s}_{2}^{\mathrm{R}}(4)\\[2pt]
		{s}_{3}^{\mathrm{R}}(1) & {s}_{3}^{\mathrm{R}}(2) & {s}_{3}^{\mathrm{R}}(3)&{s}_{3}^{\mathrm{R}}(4)\\[2pt]
		{s}_{4}^{\mathrm{R}}(1) & 
		{s}_{4}^{\mathrm{R}}(2) &{s}_{4}^{\mathrm{R}}(3)&{s}_{4}^{\mathrm{R}}(4)
	\end{bmatrix} &= \begin{bmatrix}
		\frac{1}{8} & 0&\frac{5}{8} & \frac{1}{4}\\[2pt]
		\frac{3}{5} & 0& \frac{2}{5} & 0\\[2pt]
		\frac{1}{4} & \frac{1}{2} & \frac{1}{4} & 0\\[2pt]
		\frac{2}{5} & 0 &\frac{3}{5}&0
	\end{bmatrix}.\nonumber
	\end{align}
The only time pair that satisfies the requirements stated in~\eqref{pair} is $({t},\tilde {t})=(1,3)$. As all the four jobs have fractional entries within these time intervals, we can form $\mathcal{J}'$ \purple{by choosing any $D=d^{\mathrm{max}}+1=3$ jobs out of the $4$ jobs}. Let us pick $\mathcal{J}'=\{1,2,3\}$. \purple{The next step is to form and solve \eqref{eq:xi}, which leads to, e.g., $\bm \xi = [1, -3, 2]^\top$.
One can then obtain the maximum $\Delta^\star$ satisfying~\eqref{cons1:alg2}--\eqref{cons4:alg2} for $\mathcal{J}'$ as $\Delta^\star=1/8$. Finally, we can update the schedule for the selected fractional entries by following~\eqref{cons2:5}-\eqref{cons2:6} that gives the updated schedule as:
\begin{equation}
	\mathbf{s}=\begin{bmatrix}
		\frac{1}{4} & 0&\frac{1}{2} & \frac{1}{4}\\[2pt]
		\frac{9}{40} & 0& \frac{31}{40} & 0\\[2pt]
		\frac{1}{2} & \frac{1}{2} &0& 0\\[2pt]
		\frac{2}{5} & 0 &\frac{3}{5}&0
	\end{bmatrix},\nonumber
\end{equation}
where $s_{3}(3)$ become\orange{s} integral\orangeV{.}
One can verify that the aggregate load is not changed after this adjustment. 
} 
\end{example}
\subsubsection{Analysis of the scheduling algorithm}
We next state the complexity and performance theorem for this new version of the RAR algorithm. 
\begin{lemma}[Complexity of Algorithm~\ref{alg:adjArb}]\label{lemm:1ext}
\gray{The  complexity of Algorithm~\ref{alg:adjArb} is $O\left(JT^{2}D^{3}\right)$.}
\end{lemma}

Compared with Lemma~\ref{lemm:1}, the complexity of the adjustment step is increased to accommodate time-varying loads. Nevertheless, the overall complexity is still linear in the number of jobs thus scales well for large problem instances.



\begin{theorem}[Sub-optimality bound, realistic load shapes]\label{theo:1ext}
\purple{The schedule $\mathbf {s}^\mathrm{I}$ obtained at the end of the rounding step of the RAR algorithm for realistic load shapes is feasible for problem~\eqref{opt:s}, with sub-optimality bounded as follows}:
\begin{equation}\label{bound:real}
	\Phi(\mathbf{L}^\mathrm{I}) - \Phi(\mathbf{L}^\star) \le d^{\max}T\left(T-1\right)K \max_{j\in \mathcal{J}} \big \| \mathbf{p}^{(j)} \big \|_1.
\end{equation}
\end{theorem}

Similar to the rectangular load shape case, the bound in Theorem~\ref{theo:1ext} is also independent of the number of jobs $J$.
Thus, the per-job sub-optimality approaches zero when the job set \purple{grows}. 
 \purple{The theorem is established using the same proof strategy as Theorem~\ref{theo:1} where Lemma~\ref{lemm:4} holds with an updated upper bound for Algorithm \ref{alg:adjArb}.}
%
%
%
%

	\begin{algorithm}[tp]
    \caption{Lossless adjustment, realistic loads}\label{alg:adjArb}
    \LinesNumbered
    \KwIn{$\mathbf s\gets \mathbf s^\mathrm{R}$.}

    	\While{\purple{there exists a time pair $t, \tilde t$ and a subset of jobs $\mathcal{J'}\subset \mathcal{J}$ 
    	satisfying \eqref{pair}
    	}}{
    	Find $\bm \xi\in \reals^{{D}}/\{\mathbf 0\}$ such that \eqref{eq:xi} holds.\\
    	Find $\Delta^\star$ by solving:
		\begin{subequations}\label{var:alg}
		\begin{align}
		\max_{\Delta}\quad & \Delta \\
		\mbox{s.t.} \quad 
		& s_{j_k}(t)+\xi_k{\Delta}\ge  0, \quad   j_k\in \mathcal{J'}, \label{cons1:alg2}\\
		&  s_{j_k}(t)+\xi_k{\Delta}\le  1, \quad j_k\in \mathcal{J'}, \label{cons2:alg2}\\
		& s_{j_k}(\tilde t\,)-\xi_k{\Delta}\ge  0, \quad   j_k\in \mathcal{J'}, \label{cons3:alg2}\\
		&  s_{j_k}(\tilde t\,)-\xi_k{\Delta}\le  1, \quad j_k\in \mathcal{J'}. \label{cons4:alg2}
		\end{align}
		\end{subequations}
		\\
		Update schedule $\mathbf s$:	
		\begin{align} 
		& s_{j_k}(t)\gets s_{j_k}(t)+\xi_k{\Delta^\star},\quad   j_k\in \mathcal{J'},  \label{cons2:5}\\
		& s_{j_k}(\tilde t\,)\gets s_{j_k}(\tilde t\,)-\xi_k{\Delta^\star},\quad   j_k\in \mathcal{J'}.   \label{cons2:6}
		\end{align}	

    	}
  \KwOut{$\mathbf s^\mathrm{A}\gets \mathbf s$\purple{.}}
\end{algorithm}

\section{Numerical experiments}\label{sec:ne}
In this section, we evaluate 
the performance of the proposed algorithm in the context of practical EV charging scheduling problem, where the charging power cannot be continuously adjusted. In particular, we treat the charging loads as non-preemptive, and only schedule the starting time of charging loads. We test the RAR algorithm with and without making the rectangular load shape assumption.


\subsection{Data description and preparation}
The input load parameters are derived from \purple{EV charging sessions} from Apr. 2018 to Jan. 2019 in the ACN-Data dataset \cite{lee2019acn}.
\gray{For each EV charging session, the dataset includes information on the plugged-in time (i.e., $\ta_j$), the unplugged time (i.e., $\td_j$), and the charging completion time that is between $\ta_j$ and $\td_j$.
It also provides the amount of energy delivered during the session, and the time-varying charging currents.}
\purple{
We use the length of time between plugged-in and completion time as the duration of the job. We discretize the time by 1-hour slots.}
 \purple{For realistic load shape case, we derive $\mathbf{p}^{(j)}$ from the time series of currents and the energy consumed \gray{in each charging session.} We obtain the power level for the rectangular load shape case using the average of the corresponding power profile in the realistic load shape case.}
 Over the 9557 charging sessions included, the number of EVs connected to a charger for different time of the day\footnote{We present our results in \gray{Greenwich Mean Time (GMT)}. This is consistent with the input data and circumvents the issue related with time zone change due to daylight saving.  The Pasadena local time lags GMT by 7 or 8 hours in different time of the year.} is visualized in Fig.~\ref{sol:ev}.
 

We consider a setting where the aggregator aims to match the aggregate charging load profile with that of local solar generation outputs by optimizing the following cost:
\begin{equation}
\phi_{t}\left(L(t)\right)=\left(L(t)-R(t)\right)^{2},\nonumber
\end{equation}
where $R(t)$ is the solar output within each $t \in \mathcal{T}$ derived from the solar irradiance for a sample day (see Fig.~\ref{sol:ev}) in Pasadena, California  \cite{nsrdb}.


Using these data, we create synthetic problem instances with $J \in \{100,500,2000,5000,10000\}$ charging sessions within a day. The scope of such simulated instances may represent aggregation settings arising from a single charging station to many charging stations in a city. For each $J$, we generate and solve 400 random problem instances as follows. The charging sessions are sampled with replacement from the dataset uniformly at random. The solar profile $R(t)$ is obtained by (a) sampling a realization of solar irradiance $r(t)$ from $\mathcal N(\bar r(t), \sigma^2)$ truncated to $\mathbb R_+$, where the mean irradiance $\bar r(t)$ is as described in Fig.~\ref{sol:ev} and the standard deviation $\sigma$ is set to $50\%$ of $\mathbf 1^\top \bar{\mathbf r}/T$, and (b) scaling the solar irradiance $\mathbf r$ by a constant to obtain solar output $\mathbf R$ such that $\mathbf 1^\top \mathbf R$ is $70\%$ of the aggregate EV charging loads.  

\begin{figure}[t]
	\centering
	\includegraphics[width=0.48\textwidth]{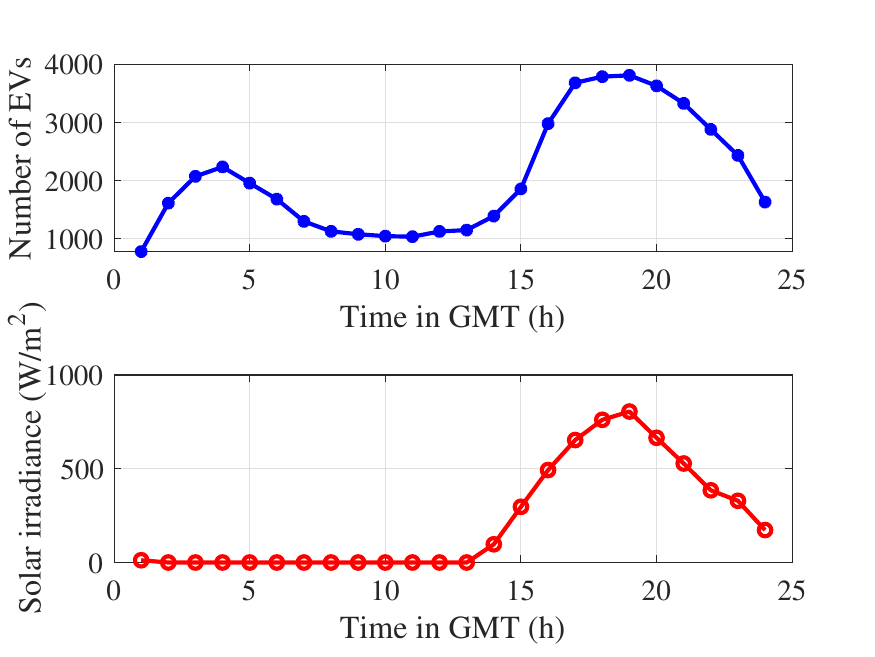}    
	\caption{Solar irradiance and number of \gray{plugged-in EVs for different time of day}}\label{sol:ev}
\end{figure}
\subsection{Deterministic scheduling}

For the generated synthetic problem instances, we evaluate the performance of the proposed algorithm and compared it with that of directly solving the MICP using Gurobi \cite{gurobi} and a greedy algorithm proposed in \cite{7066982}. All algorithms are implemented in Matlab on a Mac computer with M1 Ultra 20-core CPU and 128 GB RAM. For more consistent comparison, the convex optimization in the relaxation step of the proposed RAR algorithm is also solved with Gurobi in our experiments.


\purple{We first study the execution time of solving optimization problem~\eqref{opt:s} under the prepared settings for
both realistic and rectangular load shape cases as depicted in Fig.~\ref{Time}. The red numbers above each box shows the median of the results for that box \gray{plot.}} As shown in Fig.~\ref{Time}, although our algorithm takes more time than the MICP solver when $J$ is small, it takes an order of magnitude less time when the number of jobs is large. This superiority is more evident in the realistic load shape case as solving the MICP problem is inherently harder for solvers under this case. For example, when $J=10000$, the median run time of 
the MICP solver is 5.68 times and 12.74 times of that of the proposed algorithm for rectangular load shapes and realistic load shapes, respectively.
Similar results hold when we compare our algorithm to the greedy algorithm, with much more significant run time ratios observed for larger $J$ values. While the greedy algorithm does have a polynomial-time worst-case complexity, its process of looping over feasible schedules for each job sequentially can become very time-consuming when $J$ is large. In contrast, while the adjustment step of the RAR algorithm has a complexity linear in $J$ in the worst case, in practice, much few iterations in the adjustment steps are observed as not all the $J$ jobs will have fractional schedule in $\mathbf s^\mathrm{R}$.

%

\begin{figure}[t]
	\centering
	\begin{subfigure}{0.48\textwidth}
		\includegraphics[width=1\textwidth]{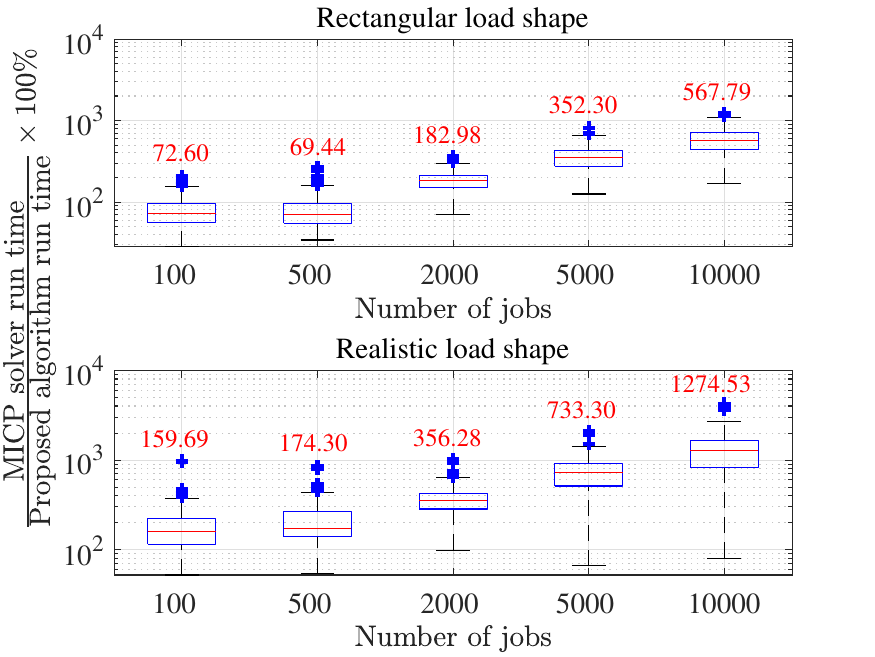}	
		\caption{Proposed algorithm v.s. MICP solver}\label{Time:milp}
	\end{subfigure}
	\begin{subfigure}{0.48\textwidth}
		\includegraphics[width=1\textwidth]{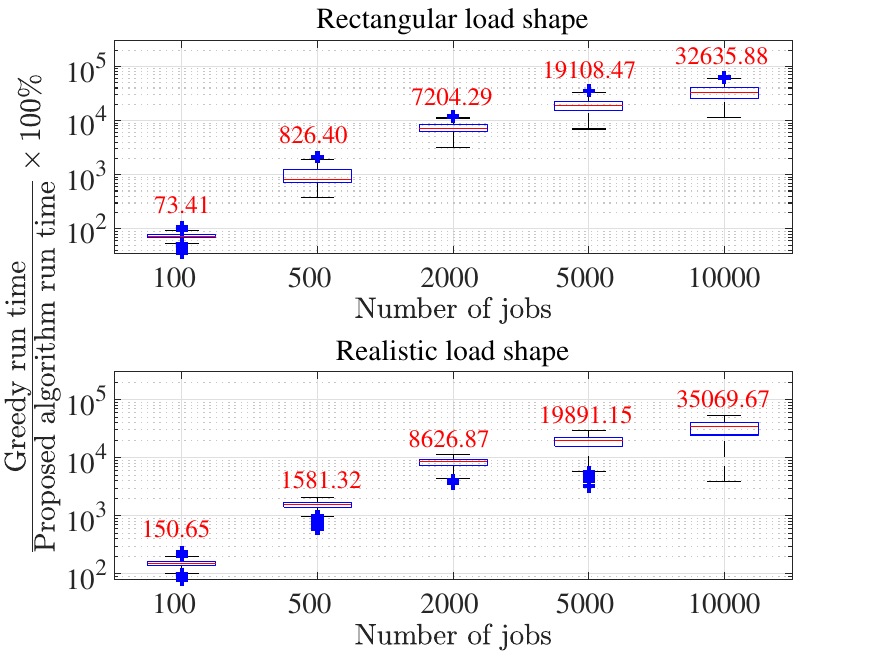}	
		\caption{Proposed algorithm v.s. Greedy}
	\end{subfigure}
	\caption{Relative computational time}\label{Time}
\end{figure}



\gray{We next compare the aggregate loads obtained by various algorithms.} As the MICP solver turned out to provide  global optimal solutions with zero sub-optimality gap in most settings (\gray{as reported by Gurobi)}, we pick its obtained aggregate load as the optimal schedule for our analysis. \gray{As an example, the aggregate loads for one experiment with $J=2000$ are depicted in Fig.~\ref{alg:sch}.} As \gray{evident} from this figure, the obtained schedule by the proposed algorithm is nearly identical to the MICP solver under both settings. \gray{The observed deviation from the MICP aggregate load with the proposed algorithm is smaller than that with greedy algorithm.} 
\begin{figure}[t]
    \centering
	\includegraphics[width=0.48\textwidth]{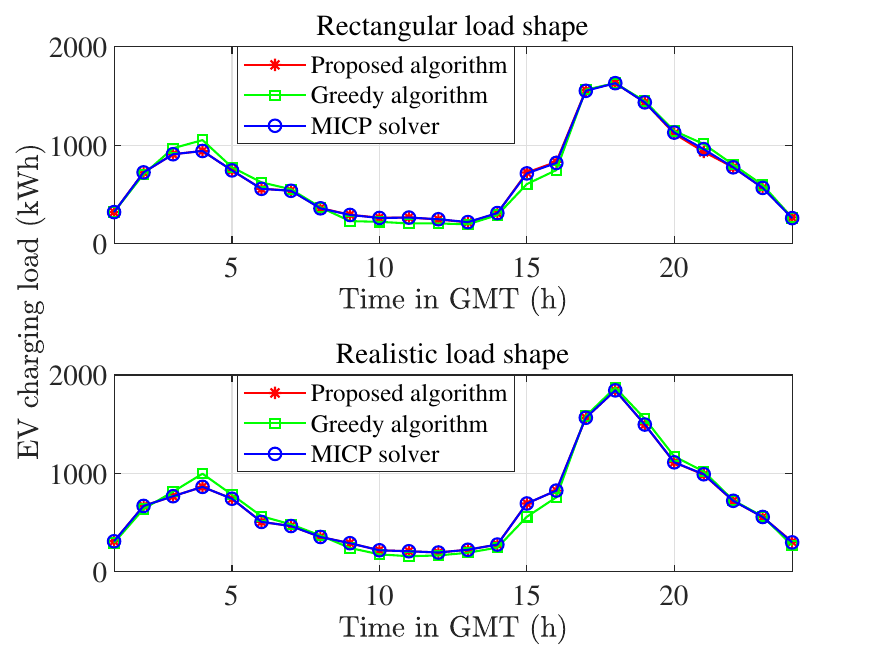}    
\caption{Scheduled aggregate load for $J=2000$}\label{alg:sch}
\end{figure}

\purple{The sub-optimality of both the proposed algorithm and the greedy algorithm is reported in Table~\ref{sub:tab}.}
For each problem instance, we evaluate the total cost obtained by the proposed algorithm,  MICP solver, and greedy algorithm denoted by  $\Phi(\mathbf{L}^\mathrm{I})$,  $\Phi(\mathbf{L}^\star)$, and  $\Phi(\mathbf{L}^\mathrm{G})$, respectively. Then, $\Phi(\mathbf{L}^\star)$ is considered as the benchmark and the sub-optimality of the aggregate load obtained by the proposed algorithm, and greedy can be evaluated by 
$
	[\Phi(\mathbf{L}^{i})-\Phi(\mathbf{L}^\star)]/{\Phi(\mathbf{L}^\star)} \times 100 \%$, for $i \in\{\mathrm I, \mathrm G\}.  
$
\begin{table}
	\scriptsize
	\centering
	\caption{Sub-optimality percentage of the algorithms}
	\label{sub:tab}
	\begin{tblr}{
			cells = {c},
			cell{1}{1} = {r=3}{},
			cell{1}{2} = {c=4}{},
			cell{1}{6} = {c=4}{},
			cell{2}{2} = {c=2}{},
			cell{2}{4} = {c=2}{},
			cell{2}{6} = {c=2}{},
			cell{2}{8} = {c=2}{},
			vlines,
			hline{1,4,9} = {-}{},
			hline{2-3} = {2-9}{},
		}
		{${J}$}     & \textbf{Proposed algorithm [\%]} &      &                    &      & \textbf{Greedy algorithm [\%]} &      &                    &      \\
		& {Rectangular}             &      & {Realistic} &      & {Rectangular}           &      & {Realistic} &      \\
		& {Mean}                             & {Std.} & {Mean}               & {Std.} & {Mean}                           & {Std.} & {Mean}               & {Std.} \\
		\hline
		$100$   & 1.03                             & 1.59 & 0.52               & 0.90 & 6.84                           & 5.04 & 11.42              & 9.94 \\
		$500$   & 0.29                             & 0.31 & 0.08               & 0.26 & 4.98                           & 2.76 & 7.33               & 3.79 \\
		$2000$  & 0.08                             & 0.07 & 0.05               & 0.17 & 4.42                           & 2.19 & 6.38               & 2.79 \\
		$5000$  & 0.03                             & 0.05 & 0.04               & 0.14 & 4.27                           & 2.00 & 6.18               & 2.61 \\
		$10000$ & 0.02                             & 0.04 & 0.04               & 0.12 & 4.19                           & 1.92 & 6.10               & 2.58 
	\end{tblr}
\end{table}

The mean and standard deviation of the \gray{percentage sub-optimality values are calculated based on the results across random problem instances}. 
Over all the $J$ values, the proposed algorithm leads to very small sub-optimality compared to the MICP solver. For $J$ large, the mean percentage sub-optimality is even no larger than $0.05\%$. Combined with the much faster run time observed previously, this suggests that the proposed algorithm can scale better than the MICP solver for large real world problem instances. 
When comparing the proposed algorithm with that of the greedy algorithm, we observe that the greedy algorithm has much higher sub-optimality percentage values.
Thus\orange{,} while both algorithms enjoy theoretical performance guarantees, the proposed algorithm empirically outperforms the greedy algorithm in our experiments.

\subsection{\textcolor{black}{Stochastic scheduling}}

\textcolor{black}{
	In this subsection, we validate the theoretical results of Section \ref{stochastic} through numerical experiments on EV scheduling with uncertainty. We consider minimizing objective function 
	$
	\mathbb{E}_{\mathbf{R}} \left[ \sum_{t \in \mathcal{T}} \left( L(t) - R(t) \right)_{+} \right],
	$
	where $R(t)$ denotes renewable generation in time $t$.
	For simplicity and to ensure Assumption \ref{assum:1} is valid, we model  \( R(t) \) for each  \( t \in \mathcal{T} \) as a stationary Gaussian random process with \( \mathcal{N}(0.66, 0.11) \), truncated to the interval \( [0\, \mathrm{MW}, 2.76 \, \mathrm{MW}] \), where the parameters are estimated using   wind data from Pasadena, California \cite{renewables_ninja}.
	For this study, we fix \( J = 3000 \).
	Note that Assumption \ref{assum:2} may be violated for this setting as $d^{\mathrm{max}} >1$.
} 

\textcolor{black}{
	To solve this stochastic scheduling problem, we adopt two distinct approaches. 
	In the \emph{first (or our proposed) approach}, leveraging the theoretical insights from~\eqref{cost:uncertain} and Theorem \ref{uncertain:result}, we replace the stochastic cost function with 
	$
	{\widetilde{\Phi}}(\mathbf{L}) := \sum_{t \in \mathcal{T}} L(t)^2,
	$
	and apply the RAR algorithm directly with this alternative cost function to compute a schedule.
	In the \emph{second (or baseline) approach}, we conduct Monte Carlo simulations to estimate the expected cost in the objective function with varying numbers of Monte Carlo samples \([25, 50, 100, 250, 500]\),  repeat each experiment 200 times, and store the resulting schedules.
	Finally, to compare the schedules from both approaches, we run a Monte Carlo simulation with sufficiently large number of samples (i.e., 20000) to evaluate their expected costs. 
}

\textcolor{black}{The cost distributions for the baseline approach are shown as box plots in Fig.~\ref{uncertain:cost}, together with the dashed green line which represents the cost for our proposed approach. As evident from the figure, it is surprising that even when Assumption \ref{assum:2} is violated, solving the stochastic scheduling problem using the \emph{baseline approach} converges to results consistent with the \emph{proposed approach}. Similar results (not reported due to page limit) are also observed when we generate the jobs ensuring \( d^{\mathrm{max}}=1 \), with which Assumption \ref{assum:2} is valid.
}

\begin{figure}[t]
	\centering
	\includegraphics[width=0.48\textwidth]{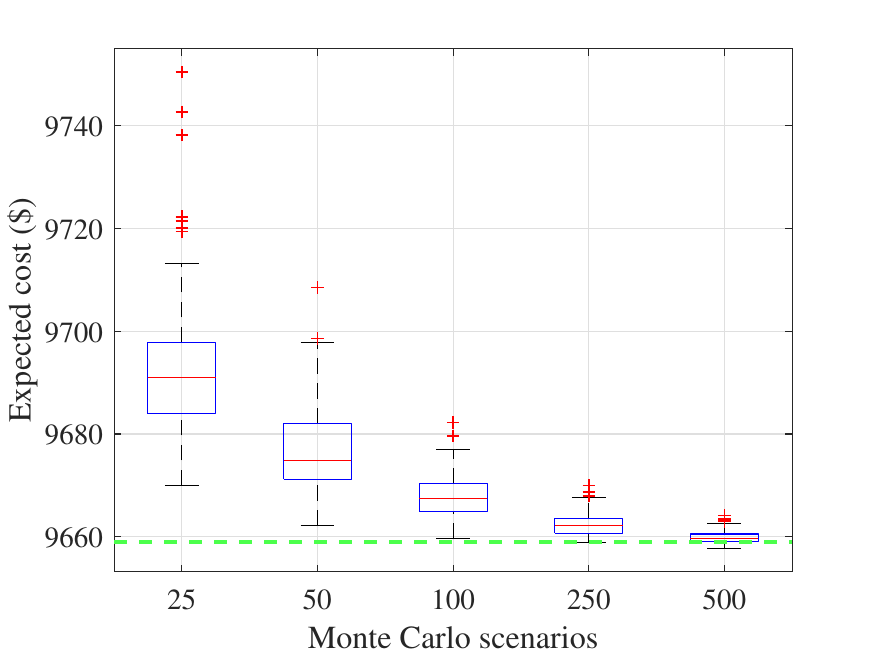}    
	\caption{\textcolor{black}{Cost distribution for stochastic scheduling: Box plots depict costs of schedules solved by using Monte Carlo  to approximate the expected cost; green dashed line corresponds to the cost obtained by the proposed modified RAR algorithm.}} \label{uncertain:cost}
\end{figure}

\section{Concluding remarks}
\gray{In this paper, we study the scheduling of non-preemptive flexible loads modeled as a mixed integer convex program. We first consider the scheduling of non-preemptive loads with rectangular load shapes and propose a polynomial-time algorithm based on a relaxation-adjustment-rounding procedure.  
We establish theoretically that the solution produced by the proposed algorithm is near-optimal for finite jobs and asymptotically optimal in a per-job cost sense when the number of jobs grows. 
The proposed RAR algorithm is then generalized to account for uncertainty and realistic time-varying load shapes, where similar performance guarantees are established.
By establishing a self-scheduling property, we  show that the RAR algorithm can be decentralized by a marginal pricing policy based on the convex relaxation of the MICP, despite the gap between the MICP and its convex relaxation. 
Overall, by designing and analyzing the novel adjustment step proposed, we are able to establish a tight connection between the MICP and its convex relaxation, offering new tools and insights for the scheduling and pricing of non-preemptive flexible loads. 
}

There are a number of interesting future directions to be explored. 
The first is the real-time version of the scheduling problem under uncertainty, where our results may serve as a bridge to enable the application of recent online convex optimization methods\orange{\cite{hazan2022introduction}}. 
The second is a more comprehensive study of the pricing problem. As discussed in Remark~\ref{rk:pricing:etc}, marginal pricing based on the convex relaxation is not the silver bullet for the general problem setting. For practically important properties like incentivizing flexibility-revealing behaviors, new pricing schemes need to be developed. 
Finally, generalizing our algorithms to more practical settings where the loads are spatially dispersed and aggregation must respect power distribution system constraints is an important area for future investigation.

%

\appendices
\section{Proofs for Section~\ref{sec:IV}} \label{appA}
\begin{IEEEproof}[Proof of Lemma~\ref{lemm:1}]
In each \orange{iteration} of the while loop, we find a cycle $\mathcal{C}$ and calculate $\Delta^\star$ and update $\mathbf{s}$. We first note that the complexity of finding a cycle is $O(T)$.
This is because an acyclic undirected multi\blue{-}graph \blue{with $T$ nodes} has at most \blue{$T-1$} edges, and therefore for any $T$ edges in $G_d(\mathbf{s})$, we can find a cycle composed of these edges using depth-first search.
Moreover, the complexity of calculating $\Delta^\star$ and updating $\mathbf{s}$ is linear with the number of edges of the cycle, which is no greater than $T$.
Thus, the complexity of each while loop \orange{iteration} is $O(T)$.

Now we analyze the number of the \orange{iterations} in the while loop.
In each \orange{iteration} of the while loop, for any $s_{j}(t)$, its value will not be changed if $s_j(t)\in \{0,1\}$.
Therefore, no edge will be created in the update.
Moreover, we note that at least one edge will be deleted in each \orange{iteration} of the while loop because there must exist $e\in \mathcal{C}$ such that either (\ref{cons:1}) or (\ref{cons:2}) is tight.
Since the total number of edges is less than $JT$ as \eqref{eq:edge:multi} specifies a spanning tree over nodes corresponding to subsets of $\mathcal T$, it follows that the adjustment algorithm terminates in $O(J T)$ \orange{iterations} \blue{across} all $d^{\max}$ multigraphs.
\end{IEEEproof}

\begin{IEEEproof}[Proof of Lemma~\ref{lemm:2}] 
In any \orange{iterations} of the while loop in Algorithm~\ref{alg:adj}, for any edge $e$ in the selected cycle, the reduction of $s_{j_e}(t_e)$ is equal to the increment of $s_{j_e}(t'_e)$ \gray{and both updated values should be within $[0, 1]$ due to~\eqref{cons:1} and~\eqref{cons:2}.}
Moreover, for any $s_{j}(t)\in \{0,1\}$, its value will not change as $(t, t', j)\not\in \mathcal{E}_{d}$ for any $t',t\in \Tcal$.
It follows that the adjustment algorithm will not make any schedule that is feasible for the convex relaxation infeasible. 
\end{IEEEproof}

\begin{IEEEproof}[Proof of Lemma~\ref{lemm:3}] 
Since $\mathbf{s}^\star$ is always a feasible schedule of the relaxed program, it is evident that $\Phi(\mathbf{L}^\mathrm{R}) \le \Phi(\mathbf{L}^\star)$.
Thus, it suffices to prove $\Phi(\mathbf{L}^\mathrm{R}) = \Phi(\mathbf{L}^\mathrm{A})$.
To do so, we show $\mathbf{L}^\mathrm{R}=\mathbf{L}^\mathrm{A}$, \blue{i.e.,} the update of $\mathbf{s}$ in Algorithm~\ref{alg:adj} will not change the aggregate load.
To this end, consider any two adjacent edges $e=(t'_e, t_e, j_e)$ and $e'=(t'_{e'}, t_{e'},j_{e'})$ in the selected cycle of any \orange{iteration} of the while loop in Algorithm~\ref{alg:adj}. 
As these two edges are adjacent, we have $t:=t_e = t'_{e'}$. 
Denote the aggregate load profile before and after the update associated with node $t$ in this \orange{iteration} of the while loop, by $\mathbf{L}$ and $\widetilde {\mathbf{L}}$, respectively. 
We also denote the corresponding schedules by $\mathbf{s}$ and $\tilde {\mathbf{s}}$, respectively. 
We have
\begin{align*}
	\widetilde {\mathbf{L}}-\mathbf{L} &= \sum_{j\in \Jcal} \mathbf{P}^{(j)} (\widetilde {\mathbf{s}}_j - \mathbf{s}_j) \\
	& = \mathbf P^{(j_e)}(\widetilde {\mathbf s}_{j_e} - \mathbf s_{j_e}) + \mathbf P^{(j_{e'})}(\widetilde {\mathbf s}_{j_{e'}} - \mathbf s_{j_{e'}})\\
	& = \mathbf P^{(j_e)} \orange{\mathbf{e}_t} \frac{\Delta^\star}{p_{j_e}} - \mathbf P^{(j_{e'})}\orange{\mathbf{e}_t} \frac{\Delta^\star}{p_{j_{e'}}} =0\orange{.}
\end{align*}

Here, the second and third identities follow from the updates associated with node $t$ in this selected cycle, and the fourth identity follows from the definition of matrices $\mathbf P^{(j)}$. 
Therefore, updates associated with node $t$ in the selected cycle \blue{do} not change the aggregate load. 
As this holds for every node in the cycle, it follows that the aggregate load is not changed after each \orange{iteration} of the while loop.
\end{IEEEproof}

\begin{IEEEproof}[Proof of Lemma~\ref{lemm:4}] 
For any $d$, the while loop terminates until there is no cycle in the graph. 
Any acyclic undirected multi\blue{-}graph \blue{with $T$ nodes} has at most \blue{$T-1$} edges. 
When the adjustment algorithm terminates, there are at most \blue{$d^{\max}(T-1)$} edges in all multigraphs, which means that there are at most \blue{$2d^{\max}(T-1) \le 2 d^{\max} T$ }fractional entries in $\mathbf s^{\mathrm{A}}$.
\end{IEEEproof}

\begin{IEEEproof}[Proof of Theorem~\ref{theo:1}] 
We start by establishing a bound for the loss introduced by the rounding step using the bound for the number of fractional entries in Lemma~\ref{lemm:4}. Note
				\begin{align*}
					\Phi(\mathbf L^\mathrm{I}) - \Phi(\mathbf L^\mathrm{A}) 
						 &\leq \sum_{t\in \mathcal{T}}  \left | \phi_t(L^\mathrm{I}(t))-\phi_t( L^\mathrm{A}(t))\right| \\
						&\le K \sum_{t\in \mathcal{T}}\left | L^\mathrm{I}(t)- L^\mathrm{A}(t)\right|,  \label{44}
				\end{align*}
				where we used the fact that $\phi_{t}$ is $K$-Lipschitz continuous for all $t\in \mathcal{T}$. We  then bound $\left | L^\mathrm{I}(t)- L^\mathrm{A}(t)\right|$ by
				\begin{align*}
				\left | L^\mathrm{I}(t)- L^\mathrm{A}(t)\right|
				&=\Big |\sum_{j\in \mathcal{J}} \sum_{t'\in \mathcal{T}}\orange{P_{t, t'}^{(j)}}\left(s_{j}^{\mathrm{I}}(t')-s_{j}^{\mathrm{A}}(t')\right)\Big| \\
				&\leq \sum_{j\in \mathcal{J}} \sum_{t'\in \mathcal{T}} \orange{P_{t, t'}^{(j)}}\left|s_{j}^{\mathrm{I}}(t')-s_{j}^{\mathrm{A}}(t')\right| \\
				&\leq \sum_{j\in \mathcal{J}} \sum_{t'\in \mathcal{T}} \orange{P_{t, t'}^{(j)}} \delta_{j,t'},
				\end{align*}
				where 
				$\delta_{j,t'}:= \mathbbm 1\{s_j^{\mathrm{I}}(t')\not=s_j^{\mathrm{A}}(t')\}$. Observe that $\delta_{j,t'}$   is an indicator having value 1 when the solution of adjustment step $s_j^\mathrm{A}(t')$ is fractional. 
It follows that
				\begin{align*}
				\Phi(\mathbf L^\mathrm{I}) - \Phi(\mathbf L^\mathrm{A})  
				& \leq  K \max_{j\in \mathcal{J}}p_jd_j \sum_{t'\in \mathcal{T}}  \sum_{j\in \mathcal{J}} \delta_{j,t'}\\
				&\leq  2d^{\max}TK \max_{j\in \mathcal{J}}p_jd_j, 
				\end{align*}
				as $\sum_{t\in \mathcal{T}} \orange{P_{t, t'}}^{(j)}\leq p_{j}d_{j}$, and $\sum_{t'\in \mathcal{T}}  \sum_{j\in \mathcal{J}} \delta_{j,t'}$ is the total number of fractional entries  in $\mathbf s^\mathrm{A}$. 
				 
Invoking Lemma~\ref{lemm:2} and Lemma~\ref{lemm:3} 
completes the proof.
\end{IEEEproof}

\begin{IEEEproof}[Proof of Lemma~\ref{adjustment:necess}] 
\blue{The lemma can be proved by constructing a problem instance where~\eqref{adjustment:necess:eq} holds. Consider the following setting:} $T\ge 2$; $J=NT$; $\phi_{t}\left(L(t)\right)= L(t)^{2}$, for all $t\in\mathcal{T}$;  $p_{j}=1$, $d_{j}=1$, and $\Tcals_j=\mathcal{T}$ for all $j\in \mathcal {J}$. It is straightforward to check that an optimal solution to the convex relaxation problem is $s_{j}^{\mathrm{R}}(t)=1/T$, \blue{for all} $t\in \mathcal{T}$ and \blue{$j\in \mathcal J$}. Then, by doing the rounding step, a possible realization of the solution, \blue{which happens with a strictly positive probability}, is $s_{j}^{\mathrm{R}}(1)=1, s_{j}^{\mathrm{R}}(t)=0$, \blue{for all  $t\in \mathcal{T}\backslash\{1\}$ and $j \in \mathcal{J}$}. This results in $L^\mathrm{I}(1)=NT$, and $L^\mathrm{I}(t)=0$,   for all $t\in \mathcal{T}\backslash\{1\}$. It is also easy to verify that $L^\star(t)=N$,  for all $t \in \mathcal{T}$. \blue{It follows that} $\Phi(\mathbf{L}^\mathrm{I})=N^{2}T^{2}$ and $\Phi(\mathbf{L}^\star)=N^{2} T$. We then have
\[
 \frac{1}{J} \left[\Phi\left(\mathbf{L}^{\mathrm{I}}\right) - \Phi\left(\mathbf{L}^\star\right) \right] = N(T-1) \ge  \frac{1}{2} NT. \IEEEQEDhereeqn
\]
\end{IEEEproof}

\section{Proofs for Section~\ref{sec:ext}} \label{appB}
\begin{IEEEproof}[Proof of Theorem~\ref{uncertain:result}]
We start by establishing a cost equivalence property for the convex relaxation of~\eqref{opt:s:u} under the imposed assumptions on the cost function. 

\begin{lemma}[Cost equivalence]\label{th:cost}
	Suppose the cost in \eqref{opt:s:u:obj} satisfies Assumptions~\ref{assum:1} and~\ref{assum:2} with some function $\psi$. Let $(\mathbf{\widetilde s}^\mathrm{R}, \mathbf{\widetilde L}^\mathrm{R})$ be a solution of the convex relaxation of~\eqref{opt:s:u} with the objective replaced by $\widetilde\Phi$. Then\orange{,} $(\mathbf{\widetilde s}^\mathrm{R}, \mathbf{\widetilde L}^\mathrm{R})$ is also a solution of the convex relaxation of~\eqref{opt:s:u}. 
\end{lemma}

\begin{IEEEproof}
As  $(\mathbf{\widetilde s}^\mathrm{R}, \mathbf{\widetilde L}^\mathrm{R})$ is a solution of the convex relaxation of~\eqref{opt:s:u} with the objective replaced by $\widetilde\Phi$, we have
	\begin{align*}
		\nabla \widetilde{\Phi}(\widetilde{\mathbf {L}}^\mathrm{R})\gray{^\top}(\mathbf L-\widetilde{\mathbf {L}}^\mathrm{R}) &= 2 (\widetilde{\mathbf {L}}^\mathrm{R})\gray{^\top} (\mathbf{L}-\widetilde{\mathbf {L}}^\mathrm{R}) \\
		&=2 (\widetilde{\mathbf {L}}^\mathrm{R})\gray{^\top} \sum_{j\in \Jcal} 	\mathbf P^{(j)} \left(\mathbf s_{j}- \widetilde{\mathbf s}_{j}^{\mathrm{R}}\right)\geq 0,
	\end{align*}
for all  $(\mathbf s, \mathbf L)$ feasible (for the convex relaxation). Since $\widetilde{\mathbf s}_{j}^\mathrm{R}$ is feasible for every $j\in \mathcal J$, the last inequality holds if and only if
\begin{equation}\label{eq:sr:ineq}
	(\widetilde{\mathbf {L}}^\mathrm{R})\gray{^\top} 	\mathbf P^{(j)} \left(\mathbf s_{j}- \widetilde{\mathbf s}_{j}^{\mathrm{R}}\right)\geq 0,
\end{equation}
for all $j \in \mathcal J$ and $\mathbf s_j$ feasible. For every $j$, \eqref{eq:sr:ineq} is equivalent to stating $\widetilde{\mathbf s}_{j}^{\mathrm{R}}$ is a solution of 
\begin{equation}\label{eq:sr:opt:form}
	\min_{\mathbf s_j \in \mathrm{conv}\, \mathcal S_j}\, (\widetilde{\mathbf {L}}^\mathrm{R})\gray{^\top} 	\mathbf P^{(j)} \mathbf s_{j},
\end{equation}
where $\mathcal S_j$ is the feasible set of $\mathbf s_j$ in \eqref{opt:s:u}. 

We claim that $\widetilde{\mathbf s}_{j}^{\mathrm{R}}$ must also be a solution of 
\begin{equation}\label{eq:sr:opt:form:Phi}
	\min_{\mathbf s_j \in \mathrm{conv}\, \mathcal S_j}\, \left(\nabla \Phi (\widetilde{\mathbf {L}}^\mathrm{R})\right)\gray{^\top} 	\mathbf P^{(j)} \mathbf s_{j}.
\end{equation}
Suppose otherwise, there must be an $\mathbf s_j'\in \mathrm{conv}\, \mathcal S_j$ such that
\[
\left(\nabla \Phi (\widetilde{\mathbf {L}}^\mathrm{R})\right)\gray{^\top} 	\mathbf P^{(j)} \widetilde{\mathbf s}_{j}^\mathrm{R} > \left(\nabla \Phi (\widetilde{\mathbf {L}}^\mathrm{R})\right)\gray{^\top} 	\mathbf P^{(j)} \mathbf s_{j}'.
\]
Without loss of generality\footnote{Indeed, for any $\mathbf s_j'$ not of the form of~\eqref{eq:pairupdate}, it is not hard to show that we can always write
$
\mathbf s_j' = \widetilde{\mathbf s}_{j}^\mathrm{R} + \sum_{k\in \mathcal K}\delta_k \orange{(\mathbf e_{t^+_k} - \mathbf e_{t^-_k})}, 
$
 where $\mathcal K$ is a finite set, $\delta_k>0$, and $t^+_k, t^-_K \in \mathcal T$ for all $k$. Furthermore, $\widetilde{\mathbf s}_{j}^\mathrm{R} +\delta_k \orange{(\mathbf e_{t^+_k} - \mathbf e_{t^-_k})}\in \mathrm{conv}\, \mathcal S_j$ for all $k\in \mathcal K$, and $\left(\nabla \Phi (\widetilde{\mathbf {L}}^\mathrm{R})\right)\gray{^\top} 	\mathbf P^{(j)} \widetilde{\mathbf s}_{j}^\mathrm{R} > \left(\nabla \Phi (\widetilde{\mathbf {L}}^\mathrm{R})\right)\gray{^\top} 	\mathbf P^{(j)} \left[\widetilde{\mathbf s}_{j}^\mathrm{R} +\delta_{k^\star} \orange{(\mathbf e_{t^+_{k^\star}} - \mathbf e_{t^-_{k^\star}})}\right]$ for some $k^\star \in \mathcal K$. We can therefore replace $\mathbf s_j'$ by $\widetilde{\mathbf s}_{j}^\mathrm{R} +\delta_{k^\star} \orange{(\mathbf e_{t^+_{k^\star}} - \mathbf e_{t^-_{k^\star}})}$ which is of the form of~\eqref{eq:pairupdate}.}, $\mathbf s_j'$ can take the form of 
\begin{equation}\label{eq:pairupdate}
	\mathbf s_j' = \widetilde{\mathbf s}_{j}^\mathrm{R} + \delta \orange{(\mathbf {e}_{t^+} - \mathbf {e}_{t^-})}, 
\end{equation}
for some $\delta >0$ and $t^+,t^- \in\mathcal T$. We thus have
\[
\delta \left(\nabla \Phi (\widetilde{\mathbf {L}}^\mathrm{R})\right)\gray{^\top} 	\mathbf P^{(j)} \orange{(\mathbf {e}_{t^+} - \mathbf {e}_{t^-})} <0.
\]
Noting that $\mathbf P^{(j)} \orange{\mathbf {e}_{t}}$ is the $t$-th column of matrix $\mathbf P^{(j)}$ and using \eqref{eq:def:Pj}, 
we have
\[
\sum_{t=t_+}^{t_++d_j-1} \psi'(\widetilde{\mathbf L}^\mathrm{R}(t)) < \sum_{t=t_-}^{t_-+d_j-1} \psi'(\widetilde{\mathbf L}^\mathrm{R}(t)).
\]
But by Assumption~\ref{assum:2}, this holds if and only if 
\[
\sum_{t=t_+}^{t_++d_j-1} \widetilde{\mathbf L}^\mathrm{R}(t) < \sum_{t=t_-}^{t_-+d_j-1} \widetilde{\mathbf L}^\mathrm{R}(t),
\]
which implies
\[
\delta \left(\widetilde{\mathbf {L}}^\mathrm{R}\right)\gray{^\top} 	\mathbf P^{(j)} \orange{(\mathbf e_{t^+} - \mathbf e_{t^-})} <0.
\]
It follows that for  $\mathbf s_j' \in \mathrm{conv}\, \mathcal S_j$ defined above, 
\[
\left(\widetilde{\mathbf {L}}^\mathrm{R}\right)\gray{^\top} 	\mathbf P^{(j)} \widetilde{\mathbf s}_{j}^\mathrm{R} > \left( \widetilde{\mathbf {L}}^\mathrm{R}\right)\gray{^\top} 	\mathbf P^{(j)} \mathbf s_{j}'.
\]
This contradicts with the fact that $\widetilde{\mathbf s}_j^\mathrm{R}$ is a solution of~\eqref{eq:sr:opt:form}. 

As $\widetilde{\mathbf s}_j^\mathrm{R}$ also optimizes \eqref{eq:sr:opt:form:Phi}, we have 
\[
\left(\nabla \Phi (\widetilde{\mathbf {L}}^\mathrm{R})\right)\gray{^\top}  	\mathbf P^{(j)} \left(\mathbf s_{j}- \widetilde{\mathbf s}_{j}^{\mathrm{R}}\right)\geq 0,
\]
for all $j \in \mathcal J$ and $\mathbf s_j$ feasible. Thus\orange{,} $(\mathbf{\widetilde s}^\mathrm{R}, \mathbf{\widetilde L}^\mathrm{R})$  also satisfies the optimality condition of the convex relaxation of \eqref{opt:s:u}, i.e., 
\[
\left(\nabla \Phi (\widetilde{\mathbf {L}}^\mathrm{R})\right)\gray{^\top}  \sum_{j \in \mathcal J}	\mathbf P^{(j)} \left(\mathbf s_{j}- \widetilde{\mathbf s}_{j}^{\mathrm{R}}\right)\geq 0,
\]
for all $j \in \mathcal J$ and $\mathbf s_j \in \mathrm{conv}\, S_j$.
\end{IEEEproof}

Equipped with this lemma, one can replace the cost function~\eqref{opt:s:u:obj} with the quadratic cost $\widetilde{\Phi}$ and then solve the resulting problem using the RAR algorithm. Applying Theorem~\ref{theo:1} gives  the sub-optimality bound which completes the proof of Theorem~\ref{uncertain:result}.
\end{IEEEproof}

\begin{remark}
Lemma~\ref{th:cost} is inspired by a similar result developed in the context of (convex) EV charging scheduling problem \cite{6313962}. In addition to re-purposing this result for a stochastic optimization application, our result is more general than that in \cite{6313962} as their problem can be shown to be a special case of the convex relaxation of our MICP when $d^{\max}=1$.
\end{remark}
\begin{IEEEproof}[Proof of Lemma~\ref{lem:price}]
We first note that
	\begin{equation}
		\pi_j^\mathrm{MP} (\mathbf s_j^\mathrm{R};\,\, \bm \lambda)=\bm \lambda^{\top} \mathbf P^{(j)}\mathbf s^\mathrm{R}_j=\sum_{t\in \mathcal{T}} \beta_j(t) s^\mathrm{R}_j(t)\gray{,}\nonumber
	\end{equation}
	where
	\begin{equation}
		\beta_j(t):=p_j \sum_{\blue{\tau\in \{t,\dots,t+d_j-1\}\cap \mathcal T}} \lambda(\tau).\nonumber
	\end{equation}
	
	Since $\mathbf s^\mathrm{R}$ is optimal for the convex relaxation, for any two periods $t_1, t_2$ such that $s^\mathrm{R}_j(t_1),s^\mathrm{R}_j(t_2)>0$, there must be $\beta_j(t_1)=\beta_j(t_2)$. 
	Otherwise, assuming $\beta_j(t_1)> \beta_j(t_2)$, we can further optimize $\mathbf s^\mathrm{R}$ by decreasing $s^\mathrm{R}_j(t_1)$ and increasing $s^\mathrm{R}_j(t_2)$.
	Indeed, define the  schedule and total aggregated load induced by this modification by $\mathbf s'$ and $\mathbf L'$, respectively. 
	By the optimality of $\mathbf s^\mathrm{R}$, we have
	$
	\nabla \Phi(\mathbf L^\mathrm{R})\gray{^\top}(\mathbf L'-\mathbf L^\mathrm{R}) \ge 0.
	$
	But
	\begin{align*}
		\nabla \Phi(\mathbf L^\mathrm{R})\gray{^\top}(\mathbf L'-\mathbf L^\mathrm{R}) &= \nabla \Phi(\mathbf L^\mathrm{R})^\top \mathbf P^{(j)}(\mathbf s_j'-\mathbf s_j^\mathrm{R}).
	\end{align*}
	Since $\bm \lambda=\nabla\Phi(\mathbf L^\mathrm{R})$, there is 
	\begin{equation}
		\nabla \Phi(\mathbf L^\mathrm{R})\gray{^\top} \mathbf P^{(j)}(\mathbf s_j'-\mathbf s_j^\mathrm{R})=\gray{\bm \beta_j^\top} (\mathbf s_j'-\mathbf s_j^\mathrm{R}).\nonumber
	\end{equation}
	By decreasing $s^\mathrm{R}_j(t_1)$ and increasing $s^\mathrm{R}_j(t_2)$ by $\epsilon>0$, we have
	\begin{equation}
		\gray{\bm\beta_j^\top} (\mathbf s_j'-\mathbf s_j^\mathrm{R})= [-\beta_j(t_1)+\beta_j(t_2)]\epsilon<0,\nonumber
	\end{equation}
	which leads to a contradiction.
	
	Given our adjustment and rounding procedure,  $s_j^{\mathrm{I}}(t)>0$ only if $s^\mathrm{R}_j(t)>0$.  It follows that
	\begin{equation}
		\sum_{t\in \mathcal{T}} \beta_j(t) s^\mathrm{R}_j(t)= \blue{\sum_{t: s^\mathrm{R}_j(t) >0} \beta_j(t) s^\mathrm{R}_j(t) = \beta_j(t^\star) \mathbf{1}^\top \mathbf s^\mathrm{R}_j=\beta_j(t^\star)},\nonumber
	\end{equation}
	where $t^\star$ denotes the period when $s^{\mathrm{I}}_j(s)=1$. 
	\blue{But  $\beta_j(t^\star) = \sum_{t\in \mathcal{T}} \beta_j(t) s^\mathrm{I}_j(t)$, thus we have} 
	\begin{equation}
		\pi_j^\mathrm{MP} (\mathbf s_j^\mathrm{R};\,\, \bm \lambda) = \pi_j^\mathrm{MP} (\mathbf s^{\mathrm{I}}_j;\,\, \bm \lambda).\nonumber \IEEEQEDhereeqn
	\end{equation}
\end{IEEEproof}
\begin{IEEEproof}[Proof of Theorem~\ref{SS}]
We start by proving a useful lemma.
\begin{lemma}
Given the marginal prices $\bm{\lambda}$, any solution of the convex relaxation of \eqref{opt:s} is also a solution of 
\begin{subequations}\label{opt:price}
	\begin{align}
		\min_{\mathbf{s}\in \reals^{J\times T},\,\, \mathbf{L} \in \reals^T} \quad &  \bm{\lambda}^\top \mathbf{L}\\
		\mbox{s.t.} \qquad\quad & \mathbf{L} = \sum_{j\in \Jcal} \mathbf{P}^{(j)} \mathbf{s}_j, \\
		& \mbox{\eqref{opt:s:c}, \eqref{opt:s:d}, and   \eqref{opt:s:bc}.}
	\end{align}
\end{subequations}
 
\end{lemma}
\begin{IEEEproof}
It is easy to verify that \blue{any $(\mathbf s, \mathbf L)$} satisfies the optimality conditions of the convex relaxation of \eqref{opt:s}, also satisfies the optimality conditions of \eqref{opt:price}. 
\end{IEEEproof}

Since $\Phi(\mathbf{L})$ and the domain of $\mathbf{L}$ are convex, for any  $\mathbf{L}'$ corresponding to a schedule $\mathbf{s}'$ feasible for the convex relaxation, by the first order optimality condition of the convex relaxation, we have
\begin{equation}
\nabla \Phi(\mathbf{L}^\mathrm{R})^\top (\mathbf{L}'-\mathbf{L}^\mathrm{R})\ge 0,\nonumber
\end{equation}
where $\mathbf{L}^\mathrm{R}$ denotes the optimal  aggregated load for the relaxed program.  
As $\bm \lambda= \nabla \Phi(\mathbf{L}^\mathrm{R})$, we have
$
\bm \lambda^\top (\mathbf{L}'-\mathbf{L}^\mathrm{R})\ge 0. 
$ 

Using this observation, for every job $j$ and admissible $\mathbf{s}_j'$, there is
\begin{equation}
\begin{aligned}
	\notag 
	&\pi_j^\mathrm{MP} (\mathbf{s}_j';\,\, \bm \lambda)-\pi_j^\mathrm{MP} (\mathbf{s}_j^\mathrm{R};\,\, \bm \lambda)
	\\=\ &\bm \lambda^{\top} \mathbf{P}^{(j)}(\mathbf{s}_j'-\mathbf{s}_j^\mathrm{R})\\=\ &\bm \lambda^{\top} \mathbf{P}^{(j)}(\mathbf{s}_j'-\mathbf{s}_j^\mathrm{R})+\sum_{i\in \mathcal{J}\setminus\{j\}}\bm \lambda^{\top} \mathbf{P}^{(i)}(\mathbf{s}_i^\mathrm{R}-\mathbf{s}_i^\mathrm{R})\\=\ &\bm \lambda^{\top} (\mathbf{L}'-\mathbf{L}^\mathrm{R})\ge 0,
\end{aligned}	
\end{equation}
\blue{where the last identity follows from the fact that $\mathbf L'=\mathbf{P}^{(j)}\mathbf{s}_j' + \sum_{i\in \mathcal{J}\setminus\{j\}} \mathbf{P}^{(i)}\mathbf{s}_i^\mathrm{R}$ indeed corresponds to a feasible schedule for the convex relaxation.}
\blue{Invoking Lemma~\ref{lem:price}, we then have, for $j\in \mathcal J$ and admissible $\mathbf{s}_j'$,}
\[
\blue{\pi_j^\mathrm{MP} (\mathbf{s}_j';\,\, \bm \lambda) \ge  \pi_j^\mathrm{MP} (\mathbf s^{\mathrm{I}}_j;\,\, \bm \lambda).} \IEEEQEDhereeqn
\]
\end{IEEEproof}

\begin{IEEEproof}[Proof of Lemma~\ref{lemm:1ext}]
\gray{The while statement (line 1) in Algorithm~\ref{alg:adjArb} can be implemented as follows. We loop over every pair of $\left(t, \tilde t\right)$ and for each pair, we initialize $J'=\emptyset$ as an empty set. We then loop over $j\in\mathcal{J}$ and check whether both $s_{j}(t)$ and $s_{j}(\tilde t)$ are fractional. If it is affirmative, we update $J'$ as $J' \cup \{j\}$. We proceed until $\left|J' \right|=D$ and apply steps in the while loop (lines 2, 3, and 4 in Algorithm~\ref{alg:adjArb}) to jobs in $J'$. After these step, at least one fractional entry in $\mathbf{s}$ is removed. We update $J'$ by removing the $j$'s such that $s_{j}(t)$ or $s_{j}(\tilde t)$ is no longer fractional. We then continue with the loop over $j \in \mathcal {J}$. Importantly, for each $(t, \tilde t)$, we can ensure no more $J'$ set of jobs satisfying the fractional condition~\eqref{pair} can be found after the updates, by looping over the set of $J$ jobs only \emph{once}. As a result, the steps in the while loop (lines 2, 3, and 4 in Algorithm~\ref{alg:adjArb}) are executed at most $O(JT^{2})$ times.

The operations in lines 2, 3, and 4 of Algorithm~\ref{alg:adjArb} have the following complexity. The complexity of finding a vector $\bm \xi$ by Gauss elimination is $O\left(D^3\right)$. Moreover, finding $\Delta^\star$ and updating $\mathbf{s}$ is 
$O\left(D\right)$. Thus, 
the overall complexity of Algorithm~\ref{alg:adjArb} is $O\left(JT^{2}D^{3}\right)$.}
\end{IEEEproof}

\begin{IEEEproof}[Proof of Theorem~\ref{theo:1ext}]
The proof follows the same major steps as that for Theorem~\ref{theo:1}. 

\subsubsection{Proof of Lemma~\ref{lemm:2} for Algorithm~\ref{alg:adjArb}}
In any cycle of the while loop in Algorithm~\ref{alg:adjArb}, for each $j_k\in \mathcal{J'}$, the increment of $s_{j_k}(t)$ is equal to the reduction of $ s_{j_k}(\tilde t)$.
Moreover, inequalities~\eqref{cons1:alg2}--\eqref{cons4:alg2} guarantee that the updated schedule is feasible. Besides, the algorithm does not change the schedule of any job $j \notin \mathcal{J'}$ or any non-fractional entries which completes the proof.
\subsubsection{Proof of Lemma~\ref{lemm:3} for Algorithm~\ref{alg:adjArb}}
\gray{We  know that  $\mathbf s^\star$ is always a feasible schedule of the relaxed problem, thus $\Phi(\mathbf L^\mathrm{R}) \le \Phi(\mathbf L^\star)$ and it remains to show  $\Phi(\mathbf L^\mathrm{R}) = \Phi(\mathbf L^\mathrm{A})$ that can be proved by showing $\mathbf L^\mathrm{R}=\mathbf L^\mathrm{A}$. Let us denote the total load before the update related to~\eqref{cons2:5}  by $\mathbf L$ and after that by $\widetilde {\mathbf L}$.} We have 
\begin{align*}
\widetilde {\mathbf L}-\mathbf L &= \sum_{j_k\in \mathcal{J'}} \mathbf P^{(j_k)} (\widetilde {\mathbf s}_{j_{k}}- \mathbf s_{j_{k}}) \\
& = \sum_{j_k\in \mathcal{J'}} \mathbf P^{(j_k)}\orange{\mathbf{e}_t}(\widetilde {s}_{j_{k}}(t)-  s_{j_{k}}(t))={\Delta^\star} \sum_{j_k\in \mathcal{J'}} \xi_k \mathbf P^{(j_k)}\orange{\mathbf{e}_t}.
\end{align*}
As $s_{j_{k}}(t)\neq 0$, $t \in \mathcal T^\mathrm{S}_{j_k}$ and in particular $t + d_{j_k} -1 \le T$ for all $j_k \in \mathcal J'$, 
we can hence write the $t$-th column of $\mathbf P^{(j_k)}$ as
\[
\mathbf P^{(j_k)}\orange{\mathbf{e}_t} = [\mathbf 0^\top, (\widehat{\mathbf{p}}^{(j_k)})^\top, \mathbf 0^\top]^\top,
\]
where the first zero vector  is of dimension $(t-1)\times 1$, and $\widehat{\mathbf{p}}^{(j_k)}$ contains the first $d_{j_k}$ elements of $\mathbf{p}^{(j_k)}$. 
In other words, $\mathbf P^{(j_k)}\orange{\mathbf{e}_t}$ is a  vector generated by padding zeros to a   possibly shortened version of  $\mathbf{p}^{(j_k)}$. In particular, there exists matrix $\mathbf H_t \in \mathbb R^{T\times d^{\max}}$ such that  $\mathbf P^{(j_k)}\orange{\mathbf{e}_t} = \mathbf H_t \mathbf{p}^{(j_k)}$. We  then have, by invoking \eqref{eq:xi},
\begin{equation}
	\widetilde {\mathbf L}-\mathbf L =\Delta^\star  \mathbf H_t \sum_{j_k\in \mathcal{J'}} \xi_k \mathbf{p}^{(j_k)} =\mathbf 0.\nonumber
\end{equation}

Thus\orange{,} the update related to~\eqref{cons2:5}  (and by the same arguments~\eqref{cons2:6}) does not change the aggregate load and therefore the cost. Since this invariance holds for every iteration of the algorithm, it must hold for the entire algorithm.

\subsubsection{Bound of the number of fractional entries} \gray{We establish a new bound on the number of fractional entries as the bound in Lemma~\ref{lemm:4} no longer applies for Algorithm~\ref{alg:adjArb}. 
Note that if a job $j\in \mathcal{J}$ has a fractional entry in $\mathbf{s}^{\mathrm{A}}_{j}$, then it has at least  two fractional entries. It can be verified that the total number of its fractional entries is less than or equal to two times the number of $(t, \tilde{t})$ pairs that both $s_{j}(t)$ and $s_{j}(\tilde{t})$ are fractional. Given that the algorithm terminates until no more set $J'$ satisfying~\eqref{pair} can be found, by pigeonhole principle, there are at most $d^{\mathrm{max}} \frac{T(T-1)}{2}$ triples of the form \orange{$(\tilde{t}, t, j)$} such that $s_j(t)$ and $s_j(\tilde{t})$ are fractional left. Thus, the upper bound for the total number of fractional entries at the end of Algorithm~\ref{alg:adjArb} is $d^{\max}T\left(T-1\right)$.} 

\subsubsection{Proof of the sub-optimality bound}
\gray{Equipped with Lemmas~\ref{lemm:2} and~\ref{lemm:3}, and the bound for the number of fractional entries, the remaining proof of this theorem is  similar to the proof of Theorem~\ref{theo:1} with two modifications: 
(a) the modified upper bound for the column sum of $P^{(j)}$, i.e., $\sum_{t\in \mathcal{T}}\orange{P^{(j)}_{t, t'}} \leq \max_{j\in \mathcal{J}} \big \| \mathbf{p}^{(j)} \big \|_1$ and (b) the new upper bound for the number fractional entries at the end of the adjustment step.  As a result,
\begin{align*}
		K\sum_{t' \in \mathcal{T}}\sum_{j\in \mathcal{J}} \delta_{j,t'} \sum_{t\in \mathcal{T}}  \orange{P_{t, t'}^{(j)}}  
	&\leq  K \max_{j\in \mathcal{J}} \big \| \mathbf{p}^{(j)} \big \|_1\sum_{t'\in \mathcal{T}}  \sum_{j\in \mathcal{J}} \delta_{j,t'}\nonumber\\
	&\leq  d^{\max}T\left(T-1\right)K \max_{j\in \mathcal{J}} \big \| \mathbf{p}^{(j)} \big \|_1 . 
\end{align*}

The remaining  proof is similar to that of Theorem~\ref{theo:1}.} 
\end{IEEEproof}

\bibliographystyle{IEEEtran}
\bibliography{bib}

\end{document}